\documentclass[twoside,12pt]{article}

\usepackage{amssymb,amsmath}

\newenvironment{proof}{\begin{trivlist}\item[]{\it
Proof.}}{\hfill$\square$\end{trivlist}}

\newtheorem{theorem}{Theorem}[section]

\newtheorem{lemma}[theorem]{Lemma}
\newtheorem{proposition}[theorem]{Proposition}

\def\mc{{\mathbb{C}}}
\def\mz{{\mathbb{Z}}}
\def\mn{{\mathbb{N}}}
\def\cuu{{\cal{U}}}
\def\coo{{\cal{O}}}
\def\caa{{\cal{A}}}
\def\ckk{{\cal{K}}}
\def\crr{{\cal{R}}}
\def\ogq{{\cal{O}}(G_q)}
\def\ugq{{\cal{U}}(G_q)}
\def\ogg{{\cal{O}}(G)}
\def\coc{^{\mathrm{coc}}}
\def\nn{{\bf{n}}}
\def\mm{{\bf{m}}}
\def\pp{{\bf{p}}}
\def\uu{{\bf{u}}}
\def\frt{{\cal{A}}(G_q)}
\def\lieg{{\mathfrak{g}}}
\def\qdet{{\cal{D}}_q}
\def\clfrt{{\cal{A}}(G)}
\def\tisi{\rho}
\def\id{{\rm{id}}}
\def\ad{{\rm{ad}}}
\def\dualad{{\rm{ad}}^{\circ}}
\def\tr{{\rm{tr}}}
\def\Tr{{\rm{Tr}}}

\def\tuq{\widetilde{U}_q} 
\def\ttt{\widetilde{T}}
\def\cartan{{\mathcal{H}}}
\def\carttau{\cartan\langle\tau\rangle} 
\def\char{{\rm{char}}}
\def\ztt{\mz[t^{\pm 1}]}
\def\res{^{{\mathrm{res}}}}
\def\ttoq{_{t\mapsto q}}
\def\ttoone{_{t\mapsto 1}}
\def\mq{{\mathbb{Q}}}
\def\rank{{\rm{rk}}}
\def\pf{{\rm{Pf}}}
\def\kkk{{\mathrm{K}}}
\def\ccc{{\mathrm{C}}}
\def\mmm{{\cal{M}}}
\def\ddd{{\cal{D}}}
\date{}

\begin{document}

\title{Representation rings of quantum groups}
\author{M Domokos ${}^a$
\thanks{This research was supported through a European
Community Marie Curie Fellowship held at the University of Edinburgh. 
Partially supported by OTKA No. F 32325, 
T 34530, 
and Leverhulme Research Interchange Grant F/00158/X.}\; 
and T H Lenagan ${}^b$
\thanks{Corresponding author.}\; 
\\ 
\\ 
{\small ${}^a$ R\'enyi Institute of Mathematics, Hungarian Academy of 
Sciences,} 
\\ {\small P.O. Box 127, 1364 Budapest, Hungary,} 
{\small E-mail: domokos@renyi.hu } 
\\ 
\\
{\small ${}^b$ 
School of Mathematics, 
University of Edinburgh,}
\\ {\small James Clerk Maxwell Building, King's Buildings,} 
{\small Mayfield Road,} 
\\ {\small Edinburgh EH9 3JZ, Scotland,}  
{\small E-mail: tom@maths.ed.ac.uk}
}

\maketitle

\begin{abstract} Generators and relations are given for the 
subalgebra of cocommutative elements in the quantized coordinate rings 
$\ogq$ of the classical groups, where $q$ is transcendental.
This is a ring theoretic formulation of the well known fact that 
the representation theory of $G_q$ is completely analogous to its 
classical counterpart. 
The subalgebras of cocommutative elements in the corresponding FRT-bialgebras 
(defined by Faddeev, Reshetikhin, and Takhtadzhyan) 
are explicitly determined, using a bialgebra embedding of 
the FRT-bialgebra into the tensor product of the quantized 
coordinate ring and the one-variable polynomial ring. 
A parallel analysis of the subalgebras of adjoint coinvariants is 
carried out as well, yielding similar results with similar proofs. 
The basic adjoint coinvariants are interpreted as 
quantum traces of representations of the corresponding quantized 
universal enveloping algebra. 
\end{abstract}



\bigskip 


\section{Introduction}\label{seq:intro} 

A good deal of classical invariant theory concerns the so called 
classical groups, their action on vectors and covectors, and their 
adjoint representation. It is therefore tempting to 
look for counterparts of this topic in the context of quantum groups,  
as is shown by various approaches in the literature. 
Our starting point here is \cite{dl1}, where two quantum versions of the 
invariant theory of the conjugation action of the general linear group
have been studied. 
Both the (right) {\it adjoint coaction}  
$\beta: f\mapsto \sum f_2\otimes S(f_1)f_3$ (given in Sweedler's notation) 
and the (right) coaction 
$\alpha:f\mapsto \sum f_2\otimes f_3S(f_1)$ 
of $\coo(GL_q(N))$, the coordinate ring of the quantum general linear group, 
on the coordinate ring of $N\times N$ quantum matrices, can be 
considered as quantum deformations of the classical conjugation action. 
In \cite{dl1}, explicit generators of the subalgebra of coinvariants were 
determined both for $\alpha$ and $\beta$ (under the assumption that 
$q$ is not a root of unity). Both algebras are $N$-variable commutative 
polynomial algebras. Note also that an element is an $\alpha$-coinvariant 
if and only if it is cocommutative. 

Some fragments of this picture had appeared in prior work already, in greater 
generality. Motivated by the theory of integrable Hamiltonian systems, 
pairwise commuting $q$-analogues of the functions $\tr(L^n)$ ($n=1,2,\dots$) 
were constructed in \cite{maillet} for algebras $\caa(R)$ 
generated by $N^2$ elements $u^i_j$, subject to the 
relations $R\uu_1\uu_2=\uu_2\uu_1R$ (see Section~\ref{sec:coc-frt} 
for explanation of this notation), where 
$R$ is an $N^2\times N^2$ matrix satisfying the Yang-Baxter equation. 
One can check that the elements constructed by Maillet are cocommutative 
in the bialgebra $\caa(R)$ (though this is not touched in \cite{maillet}). 

Another set of elements of $\caa(R)$ was constructed in 
\cite{bm}, see also \cite[Corollary 10.3.9]{majid}.   
They arise as quantum traces of powers of $u$ with respect to the 
so-called covariantized (or transmuted) product in $\caa(R)$. 
These elements are adjoint coinvariants, and pairwise commute, 
so they are also appropriate quantum analogues of the classical functions 
$\tr(L^n)$. 

The adjoint coaction is not multiplicative (neither is the version $\alpha$). 
Majid developed a theory for coquasitriangular matrix bialgebras 
$\caa(R)$ which remedies this defect. 
Namely, a new covariantized product can be introduced on 
$\caa(R)$ in a canonical way. The adjoint coinvariants become central in 
this new {\it braided matrix algebra}, known also as a 
{\it reflection equation algebra}. This process 
(called {\it transmutation} in \cite{majid}) provides a bridge between the 
results of  \cite{dl1}, and certain results on the 
reflection equation algebra. 
There is a number of papers 
dealing with the adjoint action (or coaction) on reflection 
equation algebras. 
For example, \cite{gs} and \cite{dm} (see also the references therein) 
make use of adjoint invariants 
(central elements) of the reflection equation algebra to study quantizations 
of coadjoint orbits of $SL(N)$. (Staying in the framework of quantum matrices, 
related results were obtained in \cite{dfl}.) 
See also \cite{kulish-sk} and \cite{kul-sasaki} 
for discussion of other versions of the reflection equation algebra. 
There are various versions of the Cayley-Hamilton 
theorem for quantum matrix algebras or the reflection equation algebra, 
see \cite{gps}, \cite{iop}, \cite{z}. These imply 
relations among the above mentioned adjoint coinvariants 
(respectively cocommutative elements). 

Now let us briefly describe the subject of the present paper, 
where the point of view of invariant theory is adopted, 
and we look for generators and relations for subalgebras of 
coinvariants. 
Our focus is on the matrix bialgebras $\frt$, associated with 
the classical group $G$ and the parameter $q\in\mc^{\times}$ by 
Faddeev, Reshetikhin, and Takhtadzhyan in \cite{rtf}. 
These algebras 
(called {\it FRT-bialgebras}) 
are defined in terms of generators and relations. 
They have a natural bialgebra structure, where the comultiplication 
reflects the rule for matrix multiplication. 
Following \cite{rtf}, by the coordinate ring $\ogq$ of the quantum 
group $G_q$ we mean the quotient of $\frt$ by an explicitly given ideal. 
The algebra $\frt$ 
is endowed with the adjoint coaction of $\ogq$.  
Our main result, 
Theorem~\ref{thm:coc-gen-frt} 
presents explicit generators and relations for the subalgebra 
$\frt\coc$ of cocommutative elements in $\frt$ under the assumption that 
$q$ is transcendental (the method of proof probably works when $q$ is not a 
root of unity). 
We indicate also how the same thing can be done 
for the subalgebra $\frt^{\beta}$ of adjoint 
coinvariants in $\frt$. 
This recovers the results of \cite{dl1} as the special case of 
$GL_q(N)$, $SL_q(N)$. 
For the other classical groups these results seem to be new. 
The description of $\frt\coc$ and $\frt^{\beta}$ 
is obtained from the description of the 
corresponding subalgebras in $\ogq$ (see Theorem~\ref{thm:coc-gen-ogq}), 
where the assertion is essentially a
consequence of the Peter-Weyl decomposition, due to Hayashi \cite{h}. 
Let us note however that from our point of view, the algebra $\frt$ 
is closer to the flavour of classical invariant theory (dealing with 
commutative polynomial algebras), than $\ogq$: 
it is a graded (Noetherian) algebra, having 
the same Hilbert series as its classical counterpart. 
The finite generation property of the subalgebra of coinvariants 
follows from a general Hilbert type argument, see \cite{dl2}. 

After a first draft of this paper was written, we learnt from Stephen Donkin 
that independently from us, 
strongly related results were obtained by him 
on the conjugation action of quantum groups on their 
coordinate algebra in \cite{donkin}, 
with no restriction on the deformation parameter $q$ and on the base field  
(in particular, the case when $q$ is a root of unity is covered as well). 
Moreover, his work involves the study of the structure of the coordinate ring 
of the quantum group as a module over the subalgebra of coinvariants.


\section{Cocommutative elements in $\ogq$}
\label{sec:coc-ogq} 

We work over the base field $\mc$ of complex numbers. 
Let $\ogq$ be any of the coordinate algebras of the quantum groups 
$GL_q(N)$, $SL_q(N)$, $O_q(N)$, $SO_q(N)$, $Sp_q(N)$, defined 
in sections 9.2, 9.3 of \cite{ks}, following \cite{rtf}. 
Assume that the complex parameter $q$ is not a root of unity when 
$G_q$ is $GL_q(N)$ or $SL_q(N)$, and assume that $q$ is transcendental 
in all other cases. 
We allow also the case $q=1$, when we get the 
commutative coordinate algebra $\ogg$ 
of the classical group $G$ corresponding to $G_q$. 
The assumption on $q$ guarantees that $\ogq$ is cosemisimple, 
and its corepresentation theory is completely analogous to its classical 
counterpart. The results presented in this paper depend crucially 
on the work of Hayashi \cite{h}, concerning the Peter-Weyl decomposition 
of $\ogq$. 

Recall that an element $f\in\ogq$ is cocommutative if 
$\tau\circ\Delta (f)=\Delta(f)$, 
where $\Delta:\ogq\to\ogq\otimes\ogq$ is the comultiplication, 
and $\tau$ is the flip $\tau(f\otimes g)=g\otimes f$. 
The cocommutative elements form a subalgebra $\ogq\coc$. 
We would like to point out that as an immediate corollary of the 
representation theory of $G_q$, generators and the structure of $\ogq\coc$ 
can be described explicitly. 
This is based on the following well-known statement, which is 
a reformulation of Schur's Lemma: 

\begin{lemma}\label{lem:coc-simple-coalg} 
The cocommutative elements in a simple coalgebra form a one-dimensional 
subspace. 
\end{lemma} 

\begin{proof} 
Since our base field is $\mc$, any simple coalgebra $C$ is isomorphic 
to the dual of the matrix algebra $M(N,\mc)$ for some $N$. 
The trace function on $M(N,\mc)$ fixes a vector space isomorphism 
$a\mapsto \Tr(a\cdot \_)$ between $M(N,\mc)$ and $C$. 
Under this isomorphism the center of 
$M(N,\mc)$ is mapped onto the space of cocommutative elements in $C$. 
\end{proof} 

Given a finite dimensional corepresentation 
$\varphi:V\to V\otimes\ogq$, 
write $\tr(\varphi)$ for the sum of the diagonal matrix coefficients of
$\varphi$ 
(see, for example, 1.3.2 in \cite{ks} for the notion of matrix coefficients of a 
corepresentation). 
If $\varphi$ is irreducible, then $\tr(\varphi)$ spans the 
space of cocommutative elements in the coefficient coalgebra of $\varphi$ 
by Lemma~\ref{lem:coc-simple-coalg}. 
Clearly, $\tr(\varphi\oplus\psi)=\tr(\varphi)+\tr(\psi)$, and 
$\tr(\varphi\otimes\psi)=\tr(\varphi)\cdot \tr(\psi)$, 

The isomorphism classes of irreducible corepresentations of 
$\ogq$ are parameterized by a set $P(G_q)=P(G)$. 
This set is independent of $q$, so it is the same as in the 
classical case $q=1$, when it is clearly in a natural 
bijection with the set of 
isomorphism classes of irreducible rational representations of the affine 
algebraic group $G$. 
It is a convenient tradition to represent $P(G)$ as a set of 
certain sequences of integers, 
see formulae (4.17), (6.2), and Theorem 6.4 in \cite{h}, 
or section 11.2.3 in \cite{ks} for details. 
When $G=Sp(N)$, $SL(N)$, or $GL(N)$, then it is natural to identify 
$P(G)$ with the semigroup of dominant integral weights for the corresponding 
simple Lie algebra $\lieg$, whereas when $G=SO(N)$, then $P(G)$ 
consists of those dominant integral weights for 
${\mathrm{so}}_N$, which appear as a highest weight in some 
tensor power of the vector representation of ${\mathrm{so}}_N$. 
When $G=O(N)$, then following \cite{w}, $P(G)$ is usually identified with 
the set of partitions, such that the sum of the length of the first 
two columns of their Young diagram is at most $N$. 

For $\nn\in P(G_q)$, write $\varphi_{\nn}$ 
for the corresponding irreducible corepresentation of $\ogq$. 

\begin{proposition}\label{prop:basis-coc}
The set $\{\tr(\varphi_{\nn})\mid \nn\in P(G_q)\}$ 
is a $\mc$-vector space basis of $\ogq\coc$. 
The structure constants of the algebra $\ogq\coc$ 
with respect to this basis are independent of $q$: they are the 
same as in the classical case $q=1$. 
\end{proposition} 

\begin{proof} Start with the Peter-Weyl decomposition of $\ogq$ 
due to \cite{h} 
(respectively \cite{nym} for $GL_q(N)$); see also  
11.2.3, Theorem 22 and 11.5.4, Theorem 51 in  \cite{ks}. We have  
$\ogq=\bigoplus_{\nn\in P(G_q)}C(\varphi_{\nn})$, 
where $C(\varphi_{\nn})$ is the coefficient coalgebra of 
$\varphi_{\nn}$. It follows that 
$\ogq\coc=\bigoplus_{\nn\in P(G_q)}C(\varphi_{\nn})\coc$.  
By Lemma~\ref{lem:coc-simple-coalg}, 
$C(\varphi_{\nn})\coc=\mc\tr(\varphi_{\nn})$, 
showing the first assertion. 
For the second assertion, decompose the tensor product 
$\varphi_{\nn}\otimes\varphi_{\mm}\cong 
\bigoplus_{\pp}m^{\nn,\mm}_{\pp}\varphi_{\pp}$. 
The multiplicities $m^{\nn,\mm}_{\pp}$ here are the same 
as in the classical case $q=1$, since this holds for the 
decompositions of tensor products of the corresponding 
representations of quantized universal enveloping algebras  
(see for example 7.2 in \cite{ks} or Proposition 10.1.16 in \cite{cp};  
for the case of $O_q(N)$, see  
Appendix A, Proposition~\ref{prop:tensorproduct-multiplicities} and the remark 
afterwards). 
On the other hand, they are the structure constants of $\ogq\coc$: 
we have 
$\tr(\varphi_{\nn})\cdot \tr(\varphi_{\mm})
=\sum_{\pp}m^{\nn,\mm}_{\pp}\tr(\varphi_{\pp})$. 
\end{proof} 

The following immediate corollary is a ring theoretic formulation of 
the well known fact that the representation theory of $G_q$ is essentially 
the same as the representation theory of $G$: 

\begin{proposition}\label{prop:coc-iso-classical} 
The algebra $\ogq\coc$ is isomorphic to its classical counterpart 
$\ogg\coc$, via an isomorphism mapping 
$\tr(\varphi_{\nn})\in\ogq$ to 
$\tr(\varphi_{\nn})\in\ogg$ 
for all $\nn\in P(G)$. 
\end{proposition} 

There is a natural 
right coaction of $\ogq$ on the {\it quantum exterior algebra} 
$\bigwedge(G_q)$, see sections 9.2 and 9.3 in \cite{ks}. 
The quantum exterior algebra is graded. Its degree $d$ homogeneous component 
is a subcomodule of dimension $\binom{N}{d}$, 
write $\omega_d$ for the corepresentation of $\ogq$ on this space, 
for $d=1,\ldots,N$, and set $\sigma_d=\tr(\omega_d)$. 
In the classical case $q=1$ the representation corresponding to $\omega_d$ 
is the $d$th exterior power of the defining representation of $G$. 
When $q$ is transcendental, the multiplicities of the 
irreducible summands of $\omega_d$ are the same as in the classical case 
$q=1$, since $\bigwedge(G_q)$ has the same kind of 
weight space decomposition as in the classical case.  
In particular, for $SO_q(2l)$ we have 
$\omega_l=\omega_{l,0}\oplus\omega_{l,1}$ is the direct sum of two 
non-isomorphic irreducibles; in this case set 
$\sigma_{l,0}=\tr(\omega_{l,0})$ and 
$\sigma_{l,1}=\tr(\omega_{l,1})$, so 
$\sigma_{l,0}+\sigma_{l,1}=\sigma_l$.  
Generators and relations for the commutative algebra 
$\ogq\coc$ are the following: 

\begin{theorem}\label{thm:coc-gen-ogq} 
\begin{itemize}
\item[(i)] (cf. \cite{dl2}) 
${\cal{O}}(SL_q(l+1))\coc$ is an $l$-variable commutative 
polynomial algebra generated by 
$\sigma_1,\ldots,\sigma_l$. 

\item[(ii)] ${\cal{O}}(Sp_q(2l))\coc$ is an $l$-variable 
commutative polynomial algebra 
generated by $\sigma_1,\ldots,\sigma_l$. 

\item[(iii)] ${\cal{O}}(O_q(2l+1))\coc$ is generated by 
$\sigma_1,\ldots,\sigma_l,\sigma_{2l+1}$, subject to the relation 
$\sigma_{2l+1}^2=1$. 
So it is a rank two free module generated by $1$ and $\sigma_{2l+1}$ 
over the 
$l$-variable commutative polynomial algebra 
$\mc[\sigma_1,\ldots,\sigma_l]$. 

\item[(iv)] ${\cal{O}}(SO_q(2l+1))\coc$ 
is the $l$-variable commutative polynomial algebra generated by
$\sigma_1,\ldots,\sigma_l$. 

\item[(v)] ${\cal{O}}(O_q(2l))\coc$ is generated by 
$\sigma_1,\ldots,\sigma_l,\sigma_{2l}$, subject to the relations 
$\sigma_{2l}^2=1$, $\sigma_l\sigma_{2l}=\sigma_l$. 
So it is the vector space direct sum 
$\mc[\sigma_1,\ldots,\sigma_l]
\oplus\sigma_{2l}\mc[\sigma_1,\ldots,\sigma_{l-1}]$ 
of the $l$-variable commutative polynomial algebra 
$\mc[\sigma_1,\ldots,\sigma_l]$, 
and the rank one free module generated by $\sigma_{2l}$ over the 
$(l-1)$-variable commutative polynomial algebra 
$\mc[\sigma_1,\ldots,\sigma_{l-1}]$. 

\item[(vi)] ${\cal{O}}(SO_q(2l))\coc$ is generated by 
$\sigma_1,\ldots,\sigma_{l-1},\sigma_{l,0},\sigma_{l,1}$, 
subject to the relation 
\[(\sigma_{l,0}-\sigma_{l,1})^2
=
\left(\sigma_l+2\sum_{i=0}^{l-1}\sigma_i\right)
\left(\sigma_l+2\sum_{i=0}^{l-1}(-1)^{l-i}\sigma_i\right)
,\] 
where $\sigma_l=\sigma_{l,0}+\sigma_{l,1}$. 
So ${\cal{O}}(SO_q(2l))\coc$ is a rank two free module generated by $1$ and 
$\sigma_{l,0}-\sigma_{l,1}$ 
over the $l$-variable polynomial algebra 
$\mc[\sigma_1,\ldots,\sigma_l]$. 

\item[(vii)] (cf. \cite{dl1}) 
${\cal{O}}(GL_q(N))\coc$ is the commutative 
Laurent polynomial ring generated by $\sigma_1,\ldots,\sigma_N,\sigma_N^{-1}$ 
(note that $\sigma_N$ is the quantum determinant). 

\end{itemize}
\end{theorem} 

\begin{proof} 
By Proposition~\ref{prop:coc-iso-classical} 
the result follows from its special case $q=1$. 
In the classical case the structure of $\ogg\coc$ is well known:  
it can be derived from the representation theory of $G$. 
For sake of completeness we give some references and hints in 
Appendix B. 
\end{proof} 

The quantum exterior algebra $\bigwedge(G_q)$ has a basis 
consisting of formally the same set of 
monomials as in the classical case, 
and a general monomial can be easily rewritten in terms of this 
basis, using the defining relations; 
see 9.2.1 Proposition 6, 9.3.2 Proposition 15, 9.3.4 Proposition 17 
in \cite{ks}. 
So in principle one can express the $\sigma_i$ 
for each concrete case of  
Theorem~\ref{thm:coc-gen-ogq} 
as a polynomial in the generators of $\ogq$; 
an example will be given in Section~\ref{sec:coc-frt}. 
(The cases (i) and (vii) were handled by different methods in 
\cite{dl1}, \cite{dl2}; then the $\sigma_i$ are sums of principal minors 
of the generic quantum matrix.) 
However, we do not know how to get such an expression 
for $\sigma_{l,0}$ (or $\sigma_{l,1}$) in (vi). 


\section{Cocommutative elements in the FRT-bialgebra}
\label{sec:coc-frt}

Throughout this section $G_q$ is one of $SL_q(N)$, $O_q(N)$, 
$Sp_q(N)$, and we retain the assumptions on $q$ made in 
Section~\ref{sec:coc-ogq}, so that the results of \cite{h} on the 
Peter-Weyl decomposition can be applied. 

By definition, $\ogq$ is the quotient of the so-called 
{\it FRT-bialgebra} $\frt$ modulo the ideal generated by 
$\qdet-1$, where $\qdet$ is a central group-like element, 
having degree $N$ in the case of $SL_q(N)$, and having degree $2$ 
in the cases of $O_q(N)$, $Sp_q(N)$. 
The algebra $\frt$ was defined in \cite{rtf} as the associative 
$\mc$-algebra with $N^2$ generators $u^i_j$, $(i,j=1,\dots,N)$, 
subject to the relations 
\begin{equation}\label{eq:rel-frt} 
R\uu_1\uu_2=\uu_2\uu_1R. 
\end{equation}  
Here $R$ is an $N^2\times N^2$ matrix, 
the R-matrix of the vector representation of the 
Drinfeld-Jimbo algebra $U_q(\lieg)$, 
where $\lieg$ is the simple Lie algebra corresponding to $G_q$, 
and $\uu_1=\uu\otimes I$, $\uu_2=I\otimes \uu$ are 
Kronecker products of the the $N\times N$ matrices $\uu=(u^i_j)$ and the 
identity matrix in the two possible orders. 
The relations \eqref{eq:rel-frt} are homogeneous of degree $2$ in the 
generators $u^i_j$, therefore $\frt$ is a graded algebra, 
with the generators $u^i_j$ having degree $1$. 
Moreover, $\frt$ is a bialgebra with comultiplication 
$\Delta(u^i_j)=\sum_ku^i_k\otimes u^k_j$ and counit 
$\varepsilon(u^i_j)=\delta_{i,j}$. 
Let $V$ be an $N$-dimensional $\mc$-vector space with basis 
$e_1,\ldots,e_N$. Write $\omega:V\to V\otimes\frt$ 
for the $\frt$-corepresentation given by 
$\omega(e_i)=\sum_je_j\otimes u^j_i$, and call $\omega$ the 
{\it fundamental corepresentation} of $\frt$. 
Note that the generators $u^i_j$ are nothing but 
the matrix coefficients of $\omega$ 
(with respect to the basis $e_1,\ldots,e_N$). 
It is clear then that the degree $r$ homogeneous component of $\frt$ 
is the coefficient space of the $r$th tensor power $\omega^{\otimes r}$ 
of the fundamental corepresentation. 

Write  $\pi:\frt\to\ogq$ for the natural 
surjection. 
A corepresentation $\varphi$ of $\frt$ induces the corepresentation 
$\varphi_{\ogq}=(\id\otimes\pi)\circ\varphi$ of $\ogq$. 
For $r\in\mn_0$ denote by 
$P_r(G_q)$ the subset of $P(G_q)$ consisting of the $\nn$ such that 
$(\omega^{\otimes r})_{\ogq}$,  
the $r$th tensor power of the fundamental corepresentation 
considered as a corepresentation of $\ogq$, contains a subcorepresentation 
isomorphic to $\varphi_{\nn}$; 
the explicit form of $P_r(G_q)$ can be found in 
(4.17) of \cite{h}. 
Up to isomorphism, there is a unique 
$\frt$-subcorepresentation $\varphi_{\nn,r}$ 
in $\omega^{\otimes r}$ 
with $(\varphi_{\nn,r})_{\ogq} \cong \varphi_{\nn}$. 
The coefficient space $C(\varphi_{\nn,r})$ is a simple 
subcoalgebra of the degree $r$ homogeneous component of 
$\frt$, and by 11.2.3 Theorems 21 and 22 in \cite{ks} 
we have the decomposition 
\begin{equation} \label{eq:frt-peter-weyl} 
\frt=\bigoplus_{r=0}^{\infty}\bigoplus_{\nn\in P_r(G_q)}
C(\varphi_{\nn,r}).
\end{equation}  

The polynomial ring $\mc[z]$ is a sub-bialgebra of the coordinate ring 
$\mc[z,z^{-1}]$ of the multiplicative group of $\mc$. 
The map $u^i_j\mapsto \delta_{i,j}z$ extends to a bialgebra homomorphism 
$\kappa:\frt\to \mc[z]$. This follows from the defining relations 
\eqref{eq:rel-frt}:  
specializing $\uu$ to any scalar matrix, 
$\uu_1\uu_2$ and $\uu_2\uu_1$ specialize to the 
the same scalar matrix, hence \eqref{eq:rel-frt} is fulfilled. 
Therefore there exists an algebra homomorphism $\kappa$ 
with the prescribed images of the generators. 
It is easy to check on the generators that this is a coalgebra homomorphism as 
well, moreover, that 
$\kappa$ has the following {\it centrality} property: 
\begin{equation}\label{eq:kappa-central} 
(\id\otimes\kappa)\circ\Delta_{\frt}
=\tau\circ(\kappa\otimes\id)\circ\Delta_{\frt}
\end{equation} 
where $\tau$ is the flip map $\tau(a\otimes b)= b\otimes a$. 

\begin{proposition}\label{prop:frt-injection} 
The map $\iota=(\pi\otimes\kappa)\circ\Delta_{\frt}$ is a bialgebra 
injection of $\frt$ into the tensor product bialgebra 
$\ogq\otimes\mc[z]$. The subcoalgebra $C(\varphi_{\nn,r})$ is mapped 
onto $C(\varphi_{\nn})\otimes z^r$ 
for all $r\in\mn_0$ and $\nn\in P_r(G_q)$. 
\end{proposition} 

\begin{proof} 
The map $\iota$ is defined as a composition of algebra homomorphisms, hence 
it is an algebra homomorphism. Property \eqref{eq:kappa-central} 
can be used to verify that it is a coalgebra homomorphism as well. 
The only thing left to show is that $\iota$ is injective. 
The algebra $\frt$ is graded, the generators $u^i_j$ have degree $1$. 
Similarly, the usual grading on the polynomial ring $\mc[z]$ induces a grading 
on $\ogq\otimes\mc[z]$, and the map $\iota$ is obviously homogeneous. 
Therefore the kernel of $\iota$ is spanned by homogeneous elements. 
Take a homogeneous element $f$ from $\ker(\iota)$, say of degree $r$. 
Then $\iota(f)=\pi(f)\otimes z^r$, hence $\pi(f)=0$. It follows that 
$f$ is a multiple of $\qdet-1$. 
The element $\qdet$ is not a zero-divisor 
in $\frt$ by Theorem 5.7 (1) in \cite{h}; see also   
11.2.3 Lemma 25, 
and the beginning of the proof of Theorem 22 on p.414 in \cite{ks}. 
(Note that $\frt$ is not always a domain, as we shall see later.)  
Clearly $1$ is not a zero-divisor. Therefore no non-zero multiple 
of $\qdet-1$ is homogeneous. Thus we have $f=0$. 
\end{proof} 

Write $\clfrt$ for the classical counterpart of the 
FRT-bialgebra. For $SL_q(N)$, this is just the $N^2$-variable commutative 
polynomial algebra, that we obtain when we specialize $q$ to $1$ 
in the defining relations \eqref{eq:rel-frt}. 
It is crucial to note however that in the cases of $O_q(N)$ and $Sp_q(N)$, 
the algebra $\clfrt$ is different from the $N^2$-variable 
commutative polynomial algebra  
(although specializing $q$ to $1$ in relations 
\eqref{eq:rel-frt}, 
we end up with the $N^2$-variable commutative polynomial algebra 
in these cases as well).  
To get the right definition of $\clfrt$ for $G=O(N)$ or $G=Sp(N)$, 
recall that the symmetric matrix 
$\hat R(q)=\tau\circ R$ has a spectral decomposition 
\[\hat R(q)=qP_+(q)-q^{-1}P_-(q)+\epsilon q^{\epsilon-N}P_0(q),\]  
where $\epsilon=1$ for $O_q(N)$ and $\epsilon=-1$ for $Sp_q(N)$; 
see section 9.3 in \cite{ks}. For $q$ transcendental, 
the eigenvalues $q$, $-q^{-1}$, $\epsilon q^{\epsilon-N}$ 
are pairwise different, 
therefore \eqref{eq:rel-frt} is a short expression of the equivalent set of 
relations 
\begin{equation}\label{eq:rel-frt-2}
P_+(q)\uu_1\uu_2=\uu_1\uu_2P_+(q),\ 
P_-(q)\uu_1\uu_2=\uu_1\uu_2P_-(q),\ 
P_0(q)\uu_1\uu_2=\uu_1\uu_2P_0(q).
\end{equation}  
When we specialize $q$ to 1, the eigenvalues   
$\epsilon q^{\epsilon-N}$ 
and $q$ (respectively, $-q^{-1}$) 
become equal 
in the orthogonal case (respectively, in the symplectic case), 
and that is why the relations obtained from \eqref{eq:rel-frt} are not strong 
enough. 
Instead, we can write down a third set of relations equivalent to 
\eqref{eq:rel-frt} or \eqref{eq:rel-frt-2}: 
\begin{equation}\label{eq:rel-frt-3} 
\hat R(q)\uu_1\uu_2=\uu_1\uu_2\hat R(q) \mbox{ \ and \ }
\kkk(q)\uu_1\uu_2=\uu_1\uu_2\kkk(q),  
\end{equation} 
where 
$\kkk(q)=(1+\epsilon(q-q^{-1})^{-1}(q^{N-\epsilon}-q^{\epsilon-N}))P_0(q)$. 
It is clear that \eqref{eq:rel-frt-3} is equivalent to 
\eqref{eq:rel-frt}, though \eqref{eq:rel-frt-3} 
is trivially redundant for $q\neq 1$. 
The advantage of \eqref{eq:rel-frt-3} 
compared to \eqref{eq:rel-frt-2} is that 
$\kkk(q)$ has a rather simple form.   
Write $\ccc(q)$ for the matrix of the metric defined on page 317 in 
\cite{ks}. Its non-zero entries all lie on the anti-diagonal, 
and up to sign, they are $q$-powers. 
Note that $\ccc(1)$ is the matrix of the symmetric (respectively,  
skew-symmetric) bilinear form that appears in the usual definition of 
the othogonal (respectively, symplectic) group. 
Now the entries of the $N^2\times N^2$ matrix $\kkk(q)$ are given by 
$\kkk(q)^{ji}_{mn}=\epsilon\ccc(q)^j_i\ccc(q)^m_n$, see page 318 in \cite{ks}. 
So the non-zero entries of $\kkk(q)$ are all $q$-powers up to sign. 
In particular, $K(1)$ makes sense.  
After these preparations it is natural to define $\clfrt$ as the 
algebra with generators $u^i_j$, $i,j=1,\ldots,N$, subject to the relations 
\begin{equation} \label{eq:rel-clfrt} 
\hat{R}(1)\uu_1\uu_2=\uu_1\uu_2\hat{R}(1) 
\mbox{ \ and \ }
\kkk(1)\uu_1\uu_2=\uu_1\uu_2\kkk(1).  
\end{equation} 
It is a bialgebra with comultiplication and counit given by 
the same formulae as for $\frt$. 
Specializing $q$ to $1$ in $\qdet$ we get a group-like element 
$\ddd$ of $\clfrt$. 
As we shall point out below, the quotient of $\clfrt$ modulo the ideal 
generated by $\ddd-1$ can be identified with $\ogg$, such that the images of 
the generators $u^i_j$ become the coordinate functions on $G$, 
with its usual embedding into the space $M(N,\mc)$ of $N\times N$ matrices.  

A close inspection of the proofs of the statements cited in this section 
from \cite{ks} about the coalgebra structure of $\frt$ shows that 
they remain valid for $\clfrt$. 
Indeed, the key point in the proof of \eqref{eq:frt-peter-weyl} 
is 11.2.3  Proposition 20 in \cite{ks}, which is a consequence of the quantum 
Brauer-Schur-Weyl duality, that is, that the 
commutant algebra of $\widetilde{U}_q({\mathfrak{g}})$ acting on a tensor 
power of the vector representation is generated by the `shifts' of 
$\hat{R}(q)$, see 8.6.3 Theorem 38 in \cite{ks} for a precise statement. 
Now in the classical Brauer-Schur-Weyl duality, the corresponding commutant 
algebra is generated by the shifts of $\hat{R}(1)=\tau$ and $\kkk(1)$, 
therefore the proof of 11.2.3 Proposition 20 in \cite{ks} 
works for the algebra $\clfrt$ defined in terms of $\hat{R}(1)$ and $\kkk(1)$. 
This yields a version of 11.2.3 Theorem 21 of \cite{ks} for $\clfrt$, 
and in turn the decomposition \eqref{eq:frt-peter-weyl} for $\clfrt$: 
\[ 
\clfrt=\bigoplus_{r=0}^{\infty}\bigoplus_{\nn\in P_r(G)}
C(\varphi_{\nn,r}), 
\]  
where $P_r(G)=P_r(G_q)$, since the multiplicities 
of the irreducible summands of the $r$th tensor power of the  
fundamental corepresentation of $\ogq$ are the same as for $\ogg$ 
(cf. 8.6.2 Corollary 37 (i) in \cite{ks}). 
Similarly, Proposition~\ref{prop:frt-injection} holds in the case $q=1$ as 
well. 

So we have defined $\clfrt$ as an algebra given in terms of 
generators and relations. 
The path we have followed expresses explicitly that 
$\clfrt$ is obtained as the special case $q=1$ of $\frt$. 
Moreover, we will need to compare the coalgebra structures of 
$\clfrt$ and $\frt$, and this definition makes possible a uniform approach: 
one can get the above mentioned 
statements about $\clfrt$ and the corresponding statements 
on $\frt$ with $q$ transcendental simultaneously.   
However, $\clfrt$ has a description in simple geometric terms as well. 
Namely, $\clfrt$ is the coordinate ring of the Zariski closure 
of the cone $\mc G$ of $G$, 
where by this cone we mean the image of the map 
$\mu:G\times\mc\to M(N,\mc)$, $(g,t)\mapsto tg$. 
Indeed, 
the first set of the relations \eqref{eq:rel-clfrt} 
says that the $u^i_j$ pairwise commute 
(note that $\hat{R}(1)=\tau$). By the proof of 9.3.1 Lemma 12 in \cite{ks}, 
the second set of the above relations says that 
\[\uu \ccc(1)^{-1}\uu^T \ccc(1)=\ccc(1)^{-1}\uu^T\ccc(1)\uu=
\mbox{ a scalar multiple of the identity,} \] 
where the scalar above is the quadratic group-like element $\ddd$. 
Theorems (5.2 C) and (6.3 B) of \cite{w} describe the generators 
of the vanishing ideal in the coordinate ring of $M(N,\mc)$ of 
the full orthogonal group and the symplectic group. 
This result can be paraphrased by saying that the quotient of 
$\clfrt$ modulo the ideal generated by $\ddd-1$ is indeed $\ogg$, 
as we claimed before. 
Furthermore, we obtained that the locus of solutions 
of the equations \eqref{eq:rel-clfrt} in  
$M(N,\mc)$ is 
the $N\times N$ matrix semigroup $\mmm$ consisting of the 
matrices $A$ such that  
$A\ccc(1)^{-1}A^T\ccc(1)$ and $\ccc(1)^{-1}A^T\ccc(1)A$ are equal 
scalar matrices (we allow the scalar zero). 
Clearly the subset of invertible elements in 
$\mmm$ is $\mc^{\times} G$. 
Therefore $\mmm\supseteq \overline{\mc G}$, 
there exist natural surjections 
$\pi_1:\clfrt\to\coo(\overline{\mc G})$ and 
$\pi_2:\coo(\overline{\mc G})\to\ogg$, and their composition is 
the natural surjection 
$\pi=\pi_2\circ\pi_1:\clfrt\to\ogg$. 
So, as we noted already, 
Proposition~\ref{prop:frt-injection} makes sense and is valid 
for $\clfrt$. 
It is easy to see that in this case the map $\iota$ is the composition 
$\mu^*\circ\pi_1$ of the comorphism of $\mu$ and $\pi_1$. 
Consequently, the injectivity of $\iota$ implies that 
$\pi_1$ is an isomorphism, hence 
$\mmm=\overline{\mc G}$, and 
$\clfrt$ is the coordinate ring of $\mmm$. 
(Alternatively, instead of using Proposition~\ref{prop:frt-injection}, 
it is possible to derive directly from the results of \cite{w} 
cited above that the vanishing ideal of the 
Zariski closure of $\mc G$ in $M(N,\mc)$ is generated 
by the polynomials coming from the second set of relations in 
\eqref{eq:rel-clfrt}. For sake of completeness we present this elementary 
argument in Appendix D.)  
Note that, being the coordinate ring of a linear algebraic semigroup, 
$\clfrt$ is naturally a bialgebra; 
the comultiplication and counit structures coming from this geometric 
interpretation of $\clfrt$ agree with the one specified before. 

\begin{proposition}\label{prop:frt-coc-iso-cl}
The subalgebra $\frt\coc$ of cocommutative elements in the FRT-bialgebra 
is isomorphic to its classical counterpart via an isomorphism 
mapping $\tr(\varphi_{\nn,r})\in\frt$ to $\tr(\varphi_{\nn,r})\in \clfrt$ 
for all $r\in\mn_0$, $\nn\in P_r(G)=P_r(G_q)$. 
\end{proposition} 

\begin{proof} 
By Lemma~\ref{lem:coc-simple-coalg} and \eqref{eq:frt-peter-weyl} 
we know that $\frt\coc$ has $\tr(\varphi_{\nn,r})$, 
$r\in\mn_0$, $\nn\in P_r(G)$ as a vector space basis. 
Identify $\frt$ with its image under $\iota$ from 
Proposition~\ref{prop:frt-injection}. 
Then $\frt\coc$ is identified with the subspace of 
$\ogq\coc\otimes\mc[z]$ 
spanned by the $\tr(\varphi_{\nn})\otimes z^r$ with 
$\nn\in P_r(G_q)=P_r(G)$. 
The assertion now immediately follows from 
Proposition~\ref{prop:coc-iso-classical}. 
\end{proof} 

The corepresentation 
$\omega_d$ of $\ogq$ from Section~\ref{sec:coc-ogq} is defined as 
$(\Omega_d)_{\ogq}$, where $\Omega_d$ is a natural 
right coaction of $\frt$ on the degree $d$ homogeneous component of the 
quantum exterior algebra $\bigwedge(G_q)$, for $d=1,\ldots,N$. 
Set $\tisi_d=\tr(\Omega_d)$. Then $\tisi_d$ is a cocommutative 
element in $\frt$, and $\pi(\tisi_d)=\sigma_d$. 
Another cocommutative element is $\qdet$. 
Under the bialgebra injection $\iota$, the element  
$\qdet$ is mapped to 
$1\otimes z^2$ (to $1\otimes z^N$ in the case of $SL_q(N)$), and 
$\tisi_d$ is mapped to $\sigma_d\otimes z^d$. 
The elements $\tisi_d$ can be expressed as polynomials of the 
generators $u^i_j$ in each concrete case, 
using the well known basis and the defining relations 
of $\bigwedge(G_q)$. The expression for $\qdet$ can be found in 
9.3.1 Lemma 12 of \cite{ks}. 

\bigskip
\noindent {\bf Example.} 
The quantum exterior algebra $\bigwedge(O_q(3))$ 
(we need to use the version on page 322 of \cite{ks}, 
and not the one given in \cite{rtf}) 
has three generators $y_1,y_2,y_3$, subject to the relations 
\begin{align*} 
y_1^2=y_3^2=0,\ \ \ \ y_2^2=(q^{1/2}-q^{-1/2})y_1y_3, 
\\ y_1y_2=-q^{-1} y_2y_1,\ \ y_2y_3=-q^{-1}y_3y_2,\ \ 
y_1y_3=-y_3y_1. 
\end{align*} 
For $1\leq i<j\leq 3$ we have 
$\Omega_2(y_iy_j)=\sum_{s,t=1}^3 y_sy_t\otimes u^s_iu^t_j$. 
The degree two homogeneous component of 
$\bigwedge(O_q(3))$ has the basis $y_1y_2$, $y_2y_3$, $y_1y_3$, 
and using the above relations it is easy to rewrite 
any monomial $y_sy_t$ as a linear combination of the basis elements. 
Thus one can easily get that 
\[\tisi_2=\tr(\Omega_2)
=u^1_1u^2_2-qu^2_1u^1_2
+u^2_2u^3_3-qu^3_2u^2_3
+u^1_1u^3_3-u^3_1u^1_3+(q^{1/2}-q^{-1/2})u^2_1u^2_3.\] 
An expression for the element $\qdet$ is 
\[\qdet=u^1_1u^3_3+q^{1/2}u^2_1u^2_3+qu^3_1u^1_3.\]

The explicit generators and relations for $\frt\coc$ are the following: 

\begin{theorem} \label{thm:coc-gen-frt}
\begin{itemize} 
\item[(i)] (cf. \cite{dl1}) The algebra ${\cal{A}}(SL_q(N))\coc$ 
is the $N$-variable commutative polynomial algebra generated by 
$\tisi_1,\tisi_2,\ldots,\tisi_N=\qdet$. 
In particular, its Hilbert series is 
$\prod_{i=1}^N(1-t^i)^{-1}$. 

\item[(ii)] For $Sp_q(N)$, $N=2l$, $l\in\mn$, 
the cocommutative elements ${\cal{A}}(Sp_q(N))\coc$ form 
an $(l+1)$-variable commutative polynomial algebra generated by 
$\qdet,\tisi_1,\tisi_2,\ldots,\tisi_l$. 
In particular, the Hilbert series of ${\cal{A}}(Sp_q(N))\coc$ is 
$(1-t^2)^{-1}\prod_{i=1}^l(1-t^i)^{-1}$. 

\item[(iii)] For $O_q(N)$, $N=2l\mbox{ or }2l+1$, $l\in\mn$, $N\geq 3$ 
we have that 
${\cal{A}}(O_q(N))\coc$ is the commutative algebra generated 
by $\qdet,\tisi_1,\tisi_2,\ldots,\tisi_N$, 
subject to the relations 
\begin{eqnarray*} 
\tisi_{N-i}\tisi_{N-j}=\tisi_i\tisi_j\qdet^{N-i-j}\ \ (0\leq i\leq j\leq l) \\
\tisi_i\tisi_{N-j}\qdet^{j-i}=\tisi_j\tisi_{N-i}\ \ (0\leq i<j\leq l), 
\end{eqnarray*}
where we set $\tisi_0=1$ for notational convenience. 
A $\mc$-vector space basis of ${\cal{A}}(O_q(N))\coc$ is $B(N)$, where 
\begin{eqnarray*} 
B(2l)=\{\tisi_1^{i_1}\cdots\tisi_l^{i_l}\qdet^j,\ 
\tisi_N\tisi_1^{j_1}\cdots\tisi_{l-1}^{j_{l-1}}\qdet^k,\ 
\tisi_{N-a}\qdet^b\tisi_1^{k_1}\cdots\tisi_{a-b-1}^{k_{a-b-1}}
\tisi_a^{k_a}\cdots\tisi_{l-1}^{k_{l-1}}\mid 
\\ j,k,i_s,j_s,k_s\in\mn_0,\ 0\leq b<a\leq l-1\}
\end{eqnarray*}
and 
\begin{eqnarray*}
B(2l+1)=\{\tisi_1^{i_1}\cdots\tisi_l^{i_l}\qdet^j,\ 
\tisi_N\tisi_1^{j_1}\cdots\tisi_l^{j_l}\qdet^k,\ 
\tisi_{N-a}\qdet^b 
\tisi_1^{k_1}\cdots\tisi_{a-b-1}^{k_{a-b-1}}
\tisi_a^{k_a}\cdots\tisi_l^{k_l}\mid 
\\ j,k,i_s,j_s,k_s\in\mn_0,\ 0\leq b<a\leq l\}.
\end{eqnarray*}
In particular, the Hilbert series of ${\cal{A}}(O_q(N))\coc$ is 
\[\frac{1+t^N(1-t^l)
+(1-t^2)(1-t^l)
\sum_{0\leq b<a\leq l-1}t^{N-a+2b}\prod_{k=a-b}^{a-1}(1-t^k)}
{(1-t^2)\prod_{i=1}^l(1-t^i)},\]
when $N=2l$, 
and 
\[\frac{1+t^N+(1-t^2)\sum_{0\leq b<a\leq l}
t^{N-a+2b}\prod_{k=a-b}^{a-1}(1-t^k)}
{(1-t^2)\prod_{i=1}^l(1-t^i)},\]
when $N=2l+1$. 
\end{itemize}
\end{theorem}

\begin{proof} 
By Proposition~\ref{prop:frt-coc-iso-cl} it is sufficient to prove the 
result in the classical case. 
Generators of $\clfrt\coc$ can be obtained from an old result of 
\cite{sib}. The relations among the generators can be determined using the 
classical case of 
Proposition~\ref{prop:frt-injection} and Theorem~\ref{thm:coc-gen-ogq}. 
A sketch of the details is given in Appendix C. 
\end{proof}   

The relation 
$\tisi_N^2=\qdet^N$ in $\caa(O_q(N))$  
(the special case $i=j=0$ of the first type relations in 
Theorem~\ref{thm:coc-gen-frt} (iii)) 
has already been obtained in \cite{h} and \cite{t}. 
It shows that $\caa(O_q(2l))$ is 
not a domain. 


\section{Dually paired Hopf algebras and quantum traces}
\label{sec:dually-paired} 

In this preparatory section we collect some 
standard generalities on Hopf algebras in a form that we shall need later. 

Let $\langle\cdot,\cdot\rangle:\cuu\times\coo\to\mc$ be a dual pairing of 
Hopf algebras $\cuu$ and $\coo$; 
see for example 1.2.5 in \cite{ks} for the notion of a dual pairing.
Assume that $\langle u,f\rangle=0$ for all $u$ implies $f=0$. 
Then the map $f\mapsto \langle\cdot,f\rangle$ is an injection of 
$\coo$ into the dual space $\cuu^*$ of $\cuu$. This injection identifies 
$\coo$ with a Hopf subalgebra of the finite dual $\cuu^{\circ}$ of $\cuu$; 
in the sequel we shall freely make this identification. 

Let $\varphi:V\to V\otimes\coo$, $v\mapsto \sum v_0\otimes v_1$ be a 
corepresentation of $\coo$ on $V$. (We say then that $V$ is a right 
$\coo$-comodule.) Denote by $L(V)$ the algebra of linear transformations 
on $V$. Then $\hat\varphi:\cuu\to L(V)$ defined by the formula 
$\hat\varphi(u)v:=\sum\langle u,v_1\rangle v_0$, $u\in \cuu$, $v\in V$, 
is an algebra homomorphism. Thus the corepresentation $\varphi$ on $V$ 
induces a representation $\hat\varphi$ of $\cuu$ on $V$. In other words, 
a right $\coo$-comodule $V$ automatically becomes a left $\cuu$-module, 
and the following basic properties hold:  

\begin{proposition}\label{prop:corep-rep} 
Let $\varphi:V\to V\otimes\coo$ be a corepresentation of $\coo$, and 
let $\hat\varphi$ be the corresponding representation of $\cuu$. 
\begin{itemize} 
\item[(i)] A subspace $W$ of $V$ is an $\coo$-subcomodule if and only if 
$W$ is an $\cuu$-submodule. 
\item[(ii)] An element $v\in V$ is an $\coo$-coinvariant if and only 
if $v$ is a $\cuu$-invariant. 
\item[(iii)] The coefficient space $C(\varphi)$ of $\varphi$ 
coincides with the space of matrix elements $M(\hat\varphi)$ 
of $\hat\varphi$, provided that $V$ is finite dimensional. 
\end{itemize} 
\end{proposition} 

Recall that $C(\varphi)$ is the smallest subspace $C$ in $\coo$ 
such that $\varphi(V)\subseteq V\otimes C$; it is a subcoalgebra of $\coo$. 
For a finite dimensional representation $T$ of $\cuu$ 
the space of matrix elements is
\[M(T):={\mathrm{Span}}_{\mc}\{c^{\xi}_v\mid \xi\in V^*,v\in V\}
\subset \cuu^*,\] 
where 
$V^*$ is the dual space of $V$, and for $\xi\in V^*$, $v\in V$ 
the linear function 
$c^{\xi}_v$ on $\cuu$ maps $x\in \cuu$
to $\xi(T(x)v)$. 

In the sequel we write $S$ for the antipode, and $\Delta$ 
for the comultiplication in the Hopf algebras considered. 
The {\it right adjoint coaction} $\beta:\coo\to\coo\otimes\coo$ 
is given in Sweedler's notation by 
\[\beta(f)=\sum f_2\otimes S(f_1)f_3,\] 
and the {\it right adjoint action} $\ad$ of $\cuu$ on itself is given by 
\[\ad(a)b=\sum S(a_1)ba_2,\ \quad a,b\in\cuu,\] 
see for example 1.3.4 in \cite{ks} for these definitions. 
The connection between $\ad$ and $\beta$ can be explained in terms of the 
left action $\dualad$ of $\cuu$ on its dual space $\cuu^*$, defined by the formula 
\[(\dualad(a)\xi)(b):=\xi(\ad(a)b),\ \quad a,b\in\cuu,\ \xi\in \cuu^*.\] 

\begin{proposition}\label{prop:beta-ad} 
The representation $\hat\beta$ coincides with the subrepresentation of 
$\dualad$ on the $\cuu$-invariant subspace $\coo$ of $\cuu^*$. 
\end{proposition} 

\begin{proof} 
For $a,b\in\cuu$ and $f\in\coo$ we have 
\begin{align*}
\langle b,\hat\beta(a)f\rangle 
&=\langle b,\sum\langle a,S(f_1)f_3\rangle f_2\rangle
\\ &=\sum\langle a,S(f_1)f_3\rangle\langle b,f_2\rangle
\\ &=\sum\langle a_1,S(f_1)\rangle\langle a_2,f_3\rangle
\langle b,f_2\rangle
\\ &=\sum\langle S(a_1),f_1\rangle\langle b,f_2\rangle
\langle a_2,f_3\rangle
\\ &=\sum\langle S(a_1)b,f_1\rangle\langle a_2,f_2\rangle
\\ &=\sum\langle S(a_1)ba_2,f\rangle
\\ &=\langle \ad(a)b,f\rangle.
\end{align*} 
This implies 
$\hat\beta(a)f=\ad^{\circ}(a)f$. 
\end{proof} 

Suppose that there exists an invertible element $\ckk$ 
in $\cuu$ such that 
\begin{equation}\label{eq:s-square}
S^2(a)=\ckk a\ckk^{-1}\mbox{ for all }a\in\cuu.
\end{equation}
Then for an arbitrary finite dimensional representation 
$T:\cuu\to L(V)$ we define the {\it quantum trace} of $T$ by 
\begin{equation}\label{eq:def-qtrace}
\tr_q T(a):=\Tr(T(\ckk^{-1}a)),\ \quad a\in\cuu,
\end{equation} 
where $\Tr$ is the ordinary trace function. 
So $\tr_qT$ is an element of 
$M(T)$, which is determined by the isomorphism class of $T$. 
Obviously this quantum trace depends on the choice of $\ckk$. 
It follows from \eqref{eq:s-square} and usual properties of $\Tr$ 
that 
$\ad^{\circ}(a)\tr_qT=\varepsilon(a)\tr_qT$, or, in other words, that 
$\tr_qT$ is invariant with respect to the action $\ad^{\circ}$. 

\begin{proposition}\label{prop:qtr-only-inv}
If $T:\cuu\to L(V)$ is a 
finite dimensional irreducible representation of $\cuu$, 
such that $T\otimes T^*$ and $T^*\otimes T$ are isomorphic 
representations of $\cuu$, 
then up to scalar multiple, $\tr_qT$ is the only $\ad^{\circ}$-invariant 
element in $M(T)$. 
\end{proposition} 
\begin{proof} 
We use a sequence of natural isomorphisms of $\cuu$-modules 
\begin{equation} \label{eq:nat-iso} 
L(V)\cong V\otimes V^*\cong V^*\otimes V\cong M(T). 
\end{equation} 
The first isomorphism associates with 
$v\otimes\xi\in V\otimes V^*$ the linear transformation 
$x\mapsto \xi(x)v$. This is an isomorphism of the $\cuu$-representations 
$T\otimes T^*$ and $\ad^T$, where 
\[\ad^T(a)\phi:=\sum T(a_1)\phi T(S(a_2)).\]  
By assumption, there exists a linear isomorphism 
$R_{TT^*}:V\otimes V^*\to V^*\otimes V$ 
intertwining between $T\otimes T^*$ and $T^*\otimes T$; 
this is the second isomorphism in \eqref{eq:nat-iso}. 
The third map $c:V^*\otimes V\to M(T)$ is the linear map sending 
$\xi\otimes v$ to $c^{\xi}_v$. It is surjective by the definition of 
$M(T)$.  
This map intertwines between the representations $T^*\otimes T$ and 
$\ad^{\circ}$, as one can easily check. 
(In particular, this shows that $M(T)$ is an 
$\ad^{\circ}$-invariant subspace of $\cuu^*$.) 
Since our base field $\mc$ is algebraically closed, 
the irreducibility of $T$ implies that $T(\cuu)=L(V)$, hence 
the dimension of $M(T)$ is $\dim(V)^2$. Therefore 
the surjective linear 
map $c$ goes between vector spaces of 
the same dimension. Thus $c$ must be an isomorphism. 

Since $\ckk$ is invertible, $T(\ckk^{-1})$ is non-zero, and 
so there exists a $\phi\in L(V)$ such that 
$\Tr(T(\ckk^{-1})\phi)$ is non-zero. 
Choose $a\in\cuu$ with $T(a)=\phi$. Then $\tr_qT(a)$ is non-zero, 
showing that $\tr_qT$ is a non-zero element of $M(T)$. 
Therefore 
by the $\cuu$-module isomorphisms of \eqref{eq:nat-iso}, 
it is sufficient to show that the subspace of 
$\ad^T$-invariants in $L(V)$ is one-dimensional. 
The latter statement is the assertion of Schur's lemma, because 
$\ad^T(a)\phi=\varepsilon(a)\phi$ for all $a\in\cuu$ if and only if 
$T(a)\phi=\phi T(a)$ for all $a\in\cuu$ (this equivalence can be proved by 
a straightforward modification of the well known proof of the statement 
that the center of $\cuu$ coincides with the subspace of $\ad$-invariant 
elements). 
\end{proof} 

A nice example to apply the above considerations is the case when 
$\cuu$ is {\it almost cocommutative}. This means that 
there exists an 
invertible element $\crr$ in $\cuu\otimes\cuu$ such that 
\[\tau\circ\Delta(a)=\crr\Delta(a)\crr^{-1}\mbox{ for all }a\in\cuu,\] 
where $\tau$ is the flip map. 
Set $\ckk:=\mu(\id\otimes S)(\crr^{-1})$, where 
$\mu$ is the multiplication map in $\cuu$. 
Then $\ckk$ is an invertible element of $\cuu$, with inverse 
$\mu(\id\otimes S)(\crr)$. 
Formula \eqref{eq:s-square} holds  
by Proposition 4.2.3 in \cite{cp}, and the remarks afterwards. 
Thus, using this $\ckk$, formula \eqref{eq:def-qtrace} gives an 
$\ad^{\circ}$-invariant quantum trace. 
Moreover, for an arbitrary representation $T$ of $\cuu$ 
the representations $T\otimes T^*$ and $T^*\otimes T$ are isomorphic; 
an isomorphism between them is 
$\tau\circ(T\otimes T^*)\crr$, where $\tau(v\otimes \xi)=\xi\otimes v$,   
see for example 4.2, page 119 in \cite{cp}. 
Therefore we may apply Proposition~\ref{prop:qtr-only-inv} 
to conclude that if $T$ is irreducible, then up to scalar multiple, 
$\tr_qT$ is the only $\ad^{\circ}$-invariant element in $M(T)$. 
We note that in this case 
$\tr_qT$ is the image of $\id_V\in L(V)$ under the composition of the 
isomorphisms 
\eqref{eq:nat-iso}, with the isomorphism 
$\tau\circ(T\otimes T^*)\crr$ 
being used in the middle. 


\section{Adjoint coinvariants in $\ogq$}
\label{sec:adj-coinv-ogq} 

For an arbitrary Hopf algebra $\coo$, 
the space $\coo\coc$ coincides with the space 
$\coo^{\alpha}=\{f\in \coo \mid \alpha(f)=f\otimes 1\}$ of 
$\alpha$-coinvariants, where $\alpha$ is the right coaction 
of $\coo$ on itself given in Sweedler's notation by the formula 
$\alpha:f\mapsto \sum f_2\otimes f_3 S(f_1)$, 
see \cite{dl1}. 
So in Section~\ref{sec:coc-ogq} we were dealing with $\ogq^{\alpha}$; 
a parallel analysis of the space $\ogq^{\beta}$ of $\beta$-coinvariants 
is carried out in this section, where $\beta$ is the adjoint coaction 
$\beta:f\mapsto \sum f_2\otimes S(f_1)f_3$. 
The results (and the proofs) are essentially the same as those of 
Section~\ref{sec:coc-ogq}, but the natural interpretation of them involves 
the quantized enveloping algebra $\ugq$ associated to $G_q$, 
fitting into the general framework formalized in 
Section~\ref{sec:dually-paired}. 

For $G_q=SL_q(N), Sp_q(N), SO_q(2l), SO_q(2l+1)$, 
the Hopf algebra $\ugq$ is the Drinfeld-Jimbo algebra 
$U_q({\mathrm{sl}}_N)$, $U_q({\mathrm{sp}}_N)$, 
$U_q({\mathrm{so}}_{2l})$, $U_{q^{1/2}}({\mathrm{so}}_{2l+1})$, respectively. 
The algebra $\cuu(O_q(N))$ is 
$\widetilde{U}_q({\mathrm{so}}_N)$, defined in 8.6.1 of \cite{ks}, 
following \cite{h} 
(see Appendix A of the present paper).  
The algebra $\cuu(GL_q(N))$ is $U_q({\mathrm{gl}}_N)$, defined on
page 163 of \cite{ks}. 
There is a dual pairing 
$\langle\cdot,\cdot\rangle:\ugq\times\ogq\to\mc$, 
given in 9.4 of \cite{ks}. 
We still assume that $q$ is transcendental 
(or $q$ is not a root of unity for $GL_q(N)$, $SL_q(N)$).  
Then this dual pairing is non-degenerate by \cite{h} 
(see also pages 410 and 440 in \cite{ks}). 
In particular, the map 
$f\mapsto \langle\cdot,f\rangle$ 
injects $\ogq$ into the finite dual $\ugq^{\circ}$ of 
$\ugq$.  
In the sequel we shall often consider $\ogq$ as a Hopf-subalgebra of 
$\ugq^{\circ}$ in this way. 

The representation $\hat\omega$ induced by the fundamental corepresentation 
$\omega$ is the so-called {\it vector representation} of $\ugq$. 
More generally, set 
$T_{\nn}=\hat\varphi_{\nn}$ for $\nn\in P(G_q)$.  
When $G_q=Sp_q(N)$, $SL_q(N)$, or $GL_q(N)$, then  
$\{T_{\nn}\mid \nn\in P(G_q)\}$ is a complete list of 
the isomorphism classes of the so-called 
{\it type 1 finite dimensional irreducible representations of} $\ugq$. 
When $G_q=O_q(N)$ or $SO_q(N)$, then 
$\{T_{\nn}\mid \nn\in P(G_q)\}$ is a complete list of the isomorphism 
classes of those (type 1) irreducible representations, which appear as 
a direct summand in some tensor power of the vector representation. 

Let us introduce the following ad hoc terminology. 
By the {\it basic representations of} $\ugq$ we mean 
$\hat\omega_1,\ldots,\hat\omega_l$ for $G_q=SL_q(l+1),Sp_q(2l),SO_q(2l+1)$, 
the representations 
$\hat\omega_1,\ldots,\hat\omega_N$ for $O_q(N)$, $N=2l,2l+1$, 
the representations 
$\hat\omega_1,\ldots,\hat\omega_{l-1},\hat\omega_{l,0},\hat\omega_{l,1}$ 
for $SO_q(2l)$, 
and the representations 
$\hat\omega_1,\ldots,\hat\omega_N,\hat\omega_N^*$ for $GL_q(N)$. 
 
We set 
$\ckk=K_{2\rho}\in\ugq$, where $K_{2\rho}$ 
is defined on page 164 of \cite{ks}. 
So $\rho=\sum_{i=1}^ln_i\alpha_i$ is the half-sum of positive roots, 
$\alpha_i$ are the simple roots of ${\mathfrak{g}}$, and 
$K_{2\rho}=K_1^{n_1}\cdots K_l^{n_l}$, where $K_i$ are usual generators of 
the Drinfeld-Jimbo algebra $U_q({\mathfrak{g}})$. 
For $GL_q(N)$, we set 
$\ckk= K_1^{N-1}K_2^{N-3}K_3^{N-5}\cdots K_N^{-N+1}$, 
where $K_1,\ldots,K_N$ denote the same generators of $\cuu(GL_q(N))$ 
as in 6.1.2, page 163 of \cite{ks}. 
Using 6.1.2 Proposition 6 of \cite{ks} it is easy to check that formula 
\eqref{eq:s-square} holds for $\ckk$. 
Therefore formula \eqref{eq:def-qtrace} defines an 
$\ad^{\circ}$-invariant quantum trace $\tr_qT$ for an arbitrary finite 
dimensional representation $T$ of $\ugq$. 
It is well known that for arbitrary finite dimensional representations 
$T_1,T_2$ of $\ugq$ we have $T_1\otimes T_2\cong T_2\otimes T_1$. 
Therefore by Proposition~\ref{prop:qtr-only-inv} 
we obtain that for any irreducible finite dimensional representation of 
$\ugq$, the quantum trace $\tr_qT$ spans the subspace of 
$\ad^{\circ}$-invariants in $M(T)$. 

Obviously, for finite dimensional representations $T_1, T_2$ we have 
\begin{equation}\label{eq:qtr-add}
\tr_q(T_1\oplus T_2)=\tr_qT_1+\tr_qT_2. 
\end{equation}  
Since $\ckk$ is group-like, by 7.1.6 of \cite{ks} we have 
\begin{equation}\label{eq:qtr-mult} 
\tr_q(T_1\otimes T_2)=(\tr_q T_1)\star (\tr_q T_2), 
\end{equation}
where $\star$ is the convolution multiplication in the dual of $\ugq$; 
so when the irreducible summands of 
$T_1$, $T_2$ are contained in $\{T_{\nn}\mid \nn\in P(G_q)\}$, 
then the right hand side of 
\eqref{eq:qtr-mult} is the product of 
$\tr_q T_1$ and $\tr_q T_2$ in $\ogq$. 

\begin{theorem}\label{thm:beta-coinv-ogq} 
The quantum traces $\{\tr_qT_{\nn}\mid\nn\in P(G_q)\}$ 
form a $\mc$-vector space basis of the space of $\beta$-coinvariants 
in $\ogq$. The linear map $\ogq^{\beta}\to\ogg^{\beta}$,  
$\tr_qT_{\nn}\mapsto \tr \varphi_{\nn}$,  
$\nn\in P(G)$,  
is an algebra isomorphism 
between $\ogq^{\beta}$ and its classical counterpart 
$\ogg^{\beta}=\ogg\coc$. 
As a $\mc$-algebra, $\ogq^{\beta}$ is generated by the quantum traces of 
the basic representations of $\ugq$, subject to the same 
relations as the corresponding cocommutative elements in 
Theorem~\ref{thm:coc-gen-ogq}. 
\end{theorem}

\begin{proof} Identifying $\ogq$ with a subspace of the 
dual of $\ugq$, the Peter-Weyl decomposition is written as 
$\ogq=\bigoplus_{\nn\in P(G_q)} M(T_{\nn})$. 
It is 
clearly a decomposition as a direct sum of $\beta$-subcomodules. 
Therefore we have 
$\ogq^{\beta}=\bigoplus_{\nn\in P(G_q)}M(T_{\nn})^{\beta}$, 
hence the elements $\tr_qT_{\nn}$ form a basis in $\ogq^{\beta}$ 
by Proposition~\ref{prop:qtr-only-inv}. The structure constants of 
the algebra $\ogq^{\beta}$ with respect to this basis are the multiplicities 
appearing in the tensor product decompositions 
$T_{\nn}\otimes T_{\mm}\cong 
\bigoplus_{\pp}m^{\nn,\mm}_{\pp}T_{\pp}$ 
by \eqref{eq:qtr-add} and \eqref{eq:qtr-mult}.  
Since the multiplicities $m^{\nn,\mm}_{\pp}$ are the same as in the classical 
case $q=1$ (see 7.2 of \cite{ks} or Proposition 10.1.16 in \cite{cp}, and 
Appendix A for the case of $O_q(N)$), 
we obtain the statement about the algebra 
isomorphism $\ogq^{\beta}\cong\ogg\coc$. 
Then the  
statement about the generators and relations follows from the known classical 
case (see Appendix B). 
\end{proof} 

The definition of the adjoint coaction of $\ogq$ on itself 
can be modified to make it a coaction $\beta$
of $\ogq$ on the FRT-bialgebra $\frt$ as follows: 
$\beta(f)=\sum f_2\otimes S(\pi(f_1))\pi(f_3)$. 
The results of Theorem~\ref{thm:beta-coinv-ogq} imply a description 
of $\frt^{\beta}$ both as a vector space and as an algebra with explicit 
generators and relations. This can be derived from the bialgebra embedding 
$\iota$ in Proposition~\ref{prop:frt-injection} in the same 
way as the results on $\frt\coc$. The algebra $\frt^{\beta}$ turns out to 
be isomorphic to $\frt\coc\cong\clfrt\coc=\clfrt^{\beta}$ as graded algebras 
(but $\frt^{\beta}$ and $\frt\coc$ are two different subsets of 
$\frt$ when $q\neq 1$). We omit the obvious details. 

\bigskip
\noindent{\bf Example.} Let us compute $\tr_q\hat\omega_m$ 
in the case of $GL_q(N)$. 
For subsets $I,J\subseteq \{1,\ldots,N\}$ with $|I|=|J|=m$, 
write $[I|J]$ for the corresponding quantum minor of $(u^i_j)$. 
So $[I|J]$ is the quantum determinant of the $m\times m$ quantum matrix 
$(u^i_j)^{i\in I}_{j\in J}$. 
Fix $J_0=\{1,\ldots,m\}$, and write $e_I=[J_0|I]$ for the 
quantum minors belonging to the first $m$ rows. 
Since $\Delta(e_I)=\sum_{|J|=m}e_J\otimes[J|I]$, 
the subspace in $\coo(GL_q(N))$ spanned by $\{e_J \mid m=|J|\}$ 
is a subcomodule 
with respect to the right coaction 
$\Delta$; the corresponding corepresentation is 
$\omega_m$, see 11.5.3 in \cite{ks}. 
The coefficient space $C(\omega_m)$ 
of $\omega_m$ is the subspace of $\coo(GL_q(N))$ spanned by all the 
$m\times m$ quantum minors. 
By definition of $\hat\omega_m$, for $x\in\cuu(GL_q(N))$ 
we have 
\[\hat\omega_m(x)e_I=\sum_J\langle x,[J|I]\rangle e_J.\]  
It follows from the explicit formulae giving the dual pairing in   
9.4, page 328 of \cite{ks} that 
\[\hat\omega_m(K_i)e_J= 
\begin{cases}
& q^{-1}e_J,\mbox{ if }i\in J;\\
& e_J,\mbox{ otherwise}. 
\end{cases}
\]
Consequently, we have 
\begin{align*} 
\hat\omega_m(\ckk^{-1})e_J
&=\hat\omega_m(\prod_{i=1}^NK_i^{-N-1+2i})e_J
\\&=(q^{-1})^{\sum_{i\in J}(-N-1+2i)}e_J
\\&=q^{m(N+1)}q^{-2(\sum_{i\in J}i)}e_J; 
\end{align*}
that is, the matrix of $\hat\omega_m(\ckk^{-1})$ 
with respect to the basis $\{e_J\mid m=|J|\}$ is diagonal. 
Thus 
\[\tr_q\hat\omega_m(x)
=\Tr(\hat\omega_m(\ckk^{-1})\hat\omega_m(x))
=\sum_{|J|=m}q^{(m(N+1)-2\sum_{i\in J}i)}
\langle x,[J|J]\rangle.\] 
This means that for $m=1,\ldots,N$, we have 
\[\tr_q\hat\omega_m=\sum_{|J|=m}q^{(m(N+1)-2\sum_{i\in J}i)}[J|J],\]
where the summation ranges over the $m$-element subsets $J$ of 
$\{1,\ldots,N\}$. 
Note that a scalar multiple of this element appears as the basic coinvariant 
$\tau_m$ introduced in \cite{dl1}. 
Since it is convenient to perform computations in 
$\coo(GL_q(N))$, the results of this section can be viewed 
as an explicit determination of the quantum traces of finite dimensional 
representations of type 1 of $\cuu(GL_q(N))$, as elements of 
$\coo(GL_q(N))$. 

\section{Appendix A}\label{sec:appendix} 

This appendix deals with the algebra $\cuu(O_q(2l))$ associated with 
$O_q(2l)$ by Hayashi \cite{h}. We prove the assertion on multiplicities 
of irreducibles in tensor product decompositions, used 
in the proof of Proposition~\ref{prop:basis-coc}. 

Throughout this appendix $q\in\mc^{\times}$ is not a root of unity, or $q=1$. 
Write $E_i,F_i,K_i,K_i^{-1}$ ($i=1,\ldots,l$) 
for the usual generators of the Drinfeld-Jimbo algebra 
$U_q=U_q({\mathrm{so}}_{2l})$, 
defined for example in 6.1.2 of \cite{ks}; 
when $q=1$, the algebra $U_1$ can be defined using an integral form of 
the Drinfeld-Jimbo algebra, see the proof of 
Proposition~\ref{prop:character-q-independent}.  
The universal enveloping algebra 
$U=U({\mathrm{so}}_{2l})$ 
is the homomorphic 
image of $U_1$, with kernel generated by $K_i-1$, $i=1,\ldots,l$. 
The Dynkin diagram $D_l$ of ${\mathrm{so}}_{2l}$ 
has an involutive automorphism interchanging the nodes $l-1$ and $l$. 
Denote by $\chi$ the corresponding involutive automorphism of $U_q$, 
so $E_i^{\chi}=E_{\chi(i)}$, 
$F_i^{\chi}=F_{\chi(i)}$, $K_i^{\chi}=K_{\chi(i)}$, 
where $\chi(i)=i$ for $i=1,\ldots,l-2$, 
$\chi(l-1)=\chi(l)$, and $\chi(l)=\chi(l-1)$. 
This extends to a Hopf algebra automorphism of $U_q$ 
by 6.1.6 Theorem 16 in \cite{ks}. In the case $q=1$, the automorphism 
$\chi$ of $U_1$ (see the proof of 
Proposition~\ref{prop:character-q-independent} for the definition of $\chi$ 
on $U_1$) 
induces an automorphism (denoted by $\chi$ as well) 
of the quotient $U$. (The algebra $U$ is generated by the images 
of $E_i$, $F_i$, and $\chi$ permutes them by the same rule 
as above.) 
Write $\mc[\chi]$ for the group algebra of the two-element group generated 
by $\chi$, and set $\tuq=\mc[\chi]\rtimes U_q$, 
the right crossed product algebra with commutation rule 
$\chi a\chi=a^{\chi}$, $a\in U_q$. 
Similarly, we set $\widetilde {U}=\mc[\chi]\rtimes U$. 

Let $P$ denote the weight lattice of ${\mathrm{so}}_{2l}$, and 
$P_+$ the subset of dominant integral weights. 
For $\lambda\in P_+$ denote by 
$T_{\lambda}$ the finite dimensional irreducible type 1 
representation of $U_q$ with highest weight $\lambda$. 
The type 1 finite dimensional irreducible representations of 
$\tuq$ can be determined from the corresponding list for $U_q$ 
by standard arguments, see 8.6.1 Proposition 34 in \cite{ks}. 
Namely, the automorphism $\chi$ induces an involutory action on the set of 
isomorphism classes of irreducible finite dimensional representations of 
$U_q$. Given a representation $T$ of $U_q$ define the representation 
$T^{\chi}$ by $T^{\chi}(a)=T(a^{\chi})$, $a\in U_q$. 
A dominant integral weight $\lambda$ can be represented by a sequence of 
non-negative integers 
$\lambda=(\lambda_1,\ldots,\lambda_l)$, 
where $\lambda_i=2(\lambda,\alpha_i)/(\alpha_i,\alpha_i)$, 
and $\alpha_1,\ldots,\alpha_l$ are the simple roots. 
Then $T_{\lambda}^{\chi}\cong T_{\lambda^{\chi}}$, 
where 
$\lambda^{\chi}=(\lambda_1,\ldots,\lambda_{l-2},\lambda_l,\lambda_{l-1})$. 
If $\lambda^{\chi}=\lambda$, then there are exactly two 
non-equivalent extensions of $T_{\lambda}$ to a representation of $\tuq$ 
on the same underlying space, denote them by $\ttt_{\lambda}$ 
and $\ttt_{\lambda}^{\circ}$. They are distinguished by 
$\ttt_{\lambda}(\chi)v=v$ and $\ttt_{\lambda}^{\circ}(\chi)v=-v$ 
for a highest weight vector $v$ of $T_{\lambda}$. 
If $\lambda^{\chi}\neq \lambda$, then the $U_q$-representation 
$T_{\lambda}\oplus T_{\lambda^{\chi}}$ extends to an irreducible 
representation 
$\ttt_{\lambda}$ of $\tuq$; the transformation $\ttt_{\lambda}(\chi)$ 
interchanges the underlying subspaces of 
$T_{\lambda}$ and $T_{\lambda^{\chi}}$. 
Now 
\[{\cal{T}}=\{\ttt_{\lambda},\ttt_{\lambda}^{\circ},\ttt_{\mu} \mid 
\lambda,\mu\in P_+,\ \lambda^{\chi}=\lambda,\ 
\mu^{\chi}\neq \mu\}\]  
is a complete list of isomorphism classes of type 1 finite dimensional 
irreducible representations of $\tuq$ (note that $q$ is assumed to be 
not a root of unity, or $q=1$). 

The Hopf algebra $\cuu(O_q(2l))$ was defined to be $\tuq$. 
This can be justified as follows. 
The element $\chi$ may be identified 
with a suitable reflection in the full orthogonal group $O(2l)$, 
such that the tangent map of the conjugation by 
$\chi\in O(2l)$ on the special orthogonal group $SO(2l)$,  
which is a Lie algebra automorphism of ${\mathrm{so}}_{2l}$, 
induces the automorphism $\chi$ of the universal enveloping 
algebra $U$, defined above. 
A representation $T$ of $O(2l)$ induces naturally a representation $\ttt$ of 
$\widetilde{U}$: on $U$ it is the tangent representation of 
$T$, whereas $\ttt(\chi)=T(\chi)$. Obviously, $\ttt$ determines $T$. 
Writing $P_+'$ for the subset of $P_+$ consisting of those $\lambda$ for 
which $T_{\lambda}$ is the tangent representation of a representation of 
the group $SO(2l)$ 
(note that with the notation of Section~\ref{sec:coc-ogq}, 
$P_+'$ may be naturally identified with $P(SO(2l))$), 
consider the set 
${\cal{T}}'=\{\ttt_{\lambda},\ttt_{\lambda}^{\circ},\ttt_{\mu} \mid 
\lambda,\mu\in P_+',\ \lambda^{\chi}=\lambda,\ \mu^{\chi}\neq \mu\}$ 
of $\tuq$-representations. 
In the case $q=1$, ${\cal{T}}'$ is a set  
of $\widetilde{U}$-representations (to be more precise, 
$\widetilde{U}_1$ representations factoring through 
$\widetilde{U}$), and it coincides with the set of isomorphism classes of 
$\ttt$, as $T$ ranges over the set of isomorphism classes of irreducible 
representations of $O(2l)$. 
So in the classical case $q=1$ we may think of the elements of 
${\cal{T}}'$ as representations of the full orthogonal group $O(2l)$. 

Denote by $\cartan$ the subalgebra of $U_q$ generated by 
$K_1^{\pm 1},\ldots,K_l^{\pm 1}$, 
and by $\carttau$ the subalgebra of $\tuq$ generated by 
$\chi$ over $\cartan$. 
Our next aim is to show that the structure of a finite dimensional 
$\tuq$-module is determined by its structure as an $\carttau$-module. 
Let $V$ be a type 1 finite dimensional $\tuq$-module. 
(When $q=1$, the elements $K_i$ act trivially on a type $1$ module, 
so $V$ is actually a module over $U$.) 
It has a weight space decomposition 
$V=\bigoplus_{\lambda\in P}V_{\lambda}$, 
and its $\cartan$-module structure is described by the weight multiplicities 
$(d_{\lambda}\mid\lambda\in P)$, $d_{\lambda}=\dim_{\mc}V_{\lambda}$. 
Multiplication by $\chi$ interchanges the weight spaces $V_{\lambda}$ 
and $V_{\lambda^{\chi}}$, hence $d_{\lambda}=d_{\lambda^{\chi}}$. 
For $\lambda=\lambda^{\chi}$, the subspace $V_{\lambda}$ is preserved by 
$\chi$, and $\chi$ acts as an involutory linear automorphism of $V_{\lambda}$; 
denote by $d_{\lambda}^+$ and $d_{\lambda}^-$ the multiplicity 
of $1$ and $-1$ as an eigenvalue of $\chi$ on $V_{\lambda}$, 
so $d_{\lambda}^++d_{\lambda}^-=d_{\lambda}$. 
Clearly, the $\carttau$-module structure of $V$ is determined by the 
collection of non-negative integers  
$(d_{\lambda}^+,d_{\lambda}^-,d_{\mu}\mid \lambda,\mu\in P_+,\ 
\lambda=\lambda^{\chi},\ \mu\neq\mu^{\chi})$,  
that we shall call the {\it formal character} 
$\char_{\carttau}V$ of the $\tuq$-module $V$. 

\begin{proposition}\label{prop:structure=character} 
The structure of a type 1 finite dimensional 
$\tuq$-module $V$ is determined by its 
formal character $\char_{\carttau}V$. 
\end{proposition} 

\begin{proof} By our assumption on $q$, we know that the representation $T$ 
of $\tuq$ on $V$ decomposes as a direct sum of irreducibles from 
${\cal{T}}$. One can determine this decomposition
by the following process. The formal character determines the weight 
multiplicities, hence we know how $V$ decomposes over the Drinfeld-Jimbo 
algebra $U_q$. Take a maximal weight $\lambda\in P_+$ such that $T_{\lambda}$ 
occurs with multiplicity $m>0$ in the decomposition over $U_q$. 

Case 1. If $\lambda^{\chi}\neq\lambda$, then $T_{\lambda^{\chi}}$ also occurs 
with multiplicity $m$ in the decomposition over $U_q$, and 
$T$ must contain $\ttt_{\lambda}$ as a summand with multiplicity $m$. 
Subtract $m$ times the formal character of $\ttt_{\lambda}$ from 
the formal character of $V$, and continue the same process. 

Case 2. If $\lambda^{\chi}=\lambda$, then $d_{\lambda}^++d_{\lambda}^-=m$, 
and   
$\ttt_{\lambda}$ must occur with multiplicity $d_{\lambda}^+$ in $T$, 
whereas $\ttt_{\lambda}^{\circ}$ must occur with multiplicity $d_{\lambda}^-$
in $T$. 
Subtract the formal character of these summands from $\char_{\carttau}V$, 
and continue the same process. 
\end{proof} 

For notational simplicity, set $\ttt_{\lambda}^{\circ}=\ttt_{\lambda}$, 
when $\lambda^{\chi}\neq \lambda$. 

\begin{proposition} \label{prop:character-q-independent}
The formal character of each of the $\tuq$-modules $\ttt_{\lambda}$ 
and $\ttt_{\lambda}^{\circ}$ is independent of $q$, and is the same as 
in the classical case of $\widetilde{U}$. 
\end{proposition} 

\begin{proof} 
Recall that the weight multiplicities for 
$T_{\lambda}: U_q\to {\mathrm{End}}_{\mc} V(\lambda)$ 
are independent of $q$ and are the same as in the classical case of $U$, 
see Corollary 10.1.15 in \cite{cp}. 
If $\lambda^{\chi}\neq\lambda$, 
then $\ttt_{\lambda}\big\vert_{U_q}\cong 
T_{\lambda}\oplus T_{\lambda^{\chi}}$,  
and the action of $\chi$ interchanges the weight subspaces for 
$T_{\lambda}$ and $T_{\lambda^{\chi}}$, so the assertion is obvious. 
Fix now $\lambda=\lambda^{\chi}\in P_+$, and consider 
$\ttt_{\lambda}$ (the case of $\ttt_{\lambda}^{\circ}$ is similar). 
A weight subspace $V(\lambda)_{\mu}$ for $\mu\neq\mu^{\chi}$ is interchanged 
by $\ttt_{\lambda}(\chi)$ with $V(\lambda)_{\mu^{\chi}}$, 
hence the assertion is clear for the contribution of this part in the 
formal character. 
Assume from now on that $\mu=\mu^{\chi}\in P$, 
and denote by $d^+(q)$, $d^-(q)$ the multiplicity of $+1$, $-1$  
as an eigenvalue of $\ttt_{\lambda}(\chi)$ restricted to the weight subspace 
$V(\lambda)_{\mu}$. 
What is left to show is that $d^+(q)$ and $d^-(q)$ do not depend on $q$. 

To this end we need to recall an integral form of the Drinfeld-Jimbo 
algebra. 
Let $t$ be an indeterminate, and consider the Laurent polynomial ring 
$\ztt$. Denote by $U_t=U_t^{\mq(t)}({\mathrm{so}}_{2l})$,  
the Drinfeld-Jimbo algebra over the field $\mq(t)$. 
Let $U\res$ be the $\ztt$-subalgebra of $U_t$ defined in 
9.3A of \cite{cp}. 
Write $V$ for the irreducible $U_t$-module with highest weight $\lambda$, say $v\in V$ is a fixed highest weight vector. 
Consider the $U\res$-module $V\res=U\res v$, 
and recall some of its properties from \cite[Proposition 10.1.4]{cp}. 
The module $V\res$ has a weight subspace decomposition 
$V\res=\bigoplus_{\mu\in P}V\res_{\mu}$, 
where $V\res_{\mu}$ is the intersection of $V\res$ and $V_{\mu}$. 
Moreover, each $V\res_{\mu}$ is a free $\ztt$-module. 
For $q\in\mc\setminus\{0\}$ 
define $U\res\ttoq=U\res\otimes_{\ztt}\mc$ 
and $V\res\ttoq=V\res\otimes_{\ztt}\mc$, 
where 
$\mc$ is made into a $\ztt$-module 
via the homomorphism $\ztt\to\mc$, $t\mapsto q$. 
Then $U\res\ttoq$ is the complex Drinfeld-Jimbo algebra $U_q$ 
(when $q=1$, this can be taken as the definition of $U_1$), 
and $V\res\ttoq$ is a $U_q$-module with highest weight $\lambda$ 
(for $q$ not a root of unity or $q=1$, this is the irreducible module 
associated with $\lambda$). 
Moreover, a free $\ztt$-module basis of $V\res_{\mu}$ is mapped onto a 
$\mc$-basis of the weight space $\left(V\res\ttoq\right)_{\mu}$. 

We define the automorphism $\chi$ of $U_t$ in the same way as for the complex 
Drinfeld-Jimbo algebra. Then $\chi$ permutes the $\ztt$-algebra generators 
of $U\res$, hence $\chi$ preserves $U\res$. 
The automorphism $\chi$ of $U\res$ induces an automorphism of 
$U\res\ttoq$ in an obvious manner, and the resulting automorphism clearly 
coincides with the automorphism of $U_q$ called $\chi$ already 
(when $q=1$, this can be taken as the definition of $\chi$). 
The $\mq(t)$-linear endomorphism $\ttt_{\lambda}(\chi)$ of the $U_t$-module 
$V$ preserves $V\res$. Indeed, take an element 
$a\cdot v\in V\res$, $a\in U\res $. 
Then 
\[\ttt_{\lambda}(\chi)(a\cdot v)=
(\chi\cdot a)\cdot v=(a^{\chi}\cdot \chi)\cdot v
=a^{\chi}\cdot (\ttt_{\lambda}(\chi)(v))
=\pm a^{\chi}\cdot v\in U\res v,\]
since $a^{\chi}\in U\res$. 
So $\ttt_{\lambda}(\chi)$ restricts to an automorphism of 
the free $\ztt$-module $V\res_{\mu}$ (recall that $\mu^{\chi}=\mu$). 
This free $\ztt$-module automorphism is represented by a square matrix 
$B$ with entries from $\ztt$, such that $B^2=I$, the identity matrix. 
Denote by $B\ttoone$ the integer matrix obtained from $B$ by specializing $t$ 
to $1$. Then $B\ttoone$ is a matrix which represents $\chi$, 
acting on $\left(V\res\ttoone\right)_{\mu}$ 
via the classical irreducible representation $\ttt_{\lambda}$ 
of $\widetilde{U}$. 
Clearly we have $\rank_{\mq(t)}(B-I)\geq \rank_{\mc}(B\ttoone-I)$ 
and $\rank_{\mq(t)}(B+I)\geq\rank_{\mc}(B\ttoone+I)$, implying 
$d^+(t)\leq d^+(1)$ and $d^-(t)\leq d^-(1)$. 
On the other hand, 
$d^+(t)+d^-(t)=
\dim_{\mq(t)}V_{\mu}
=\dim_{\mc}\left(V\res\ttoone\right)_{\mu}
=d^+(1)+d^-(1)$, 
so we have equality in both of the above inequalities. 
Similarly, for $q\in\mc\setminus\{0\}$ 
write $B\ttoq$ for the complex matrix obtained from $B$ by specializing 
$t$ to $q$. Then $B\ttoq$ represents $\chi$, acting on 
$\left(V\res\ttoq\right)_{\mu}$ 
via the representation $\ttt_{\lambda}$ of $\tuq$. 
The obvious inequalities 
$\rank_{\mq(t)}(B-I)\geq \rank_{\mc}(B\ttoq-I)$ and 
$\rank_{\mq(t)}(B+I)\geq\rank_{\mc}(B\ttoq+I)$ imply  
$d^+(t)\leq d^+(q)$ and $d^-(t)\leq d^-(q)$. 
On the other hand, 
$d^+(t)+d^-(t)
=\dim_{\mq(t)}V_{\mu}
=\dim_{\mc}\left(V\res\ttoq\right)_{\mu}
=d^+(q)+d^-(q)$. 
Hence we get 
$d^+(q)=d^+(t)=d^+(1)$ 
and $d^-(q)=d^-(t)=d^-(1)$. 
\end{proof} 

\begin{proposition} \label{prop:tensorproduct-multiplicities} 
Assume that $q\in\mc\setminus\{0\}$ is not a root of unity or $q=1$. 
Then the tensor product of any pair of representations 
$T_1,T_2\in{\cal{T}}$ decomposes as 
\[T_1\otimes T_2\cong\bigoplus_{T\in{\cal{T}}}m_TT,\] 
and the multiplicities $m_T$ here are  independent of $q$ 
(they are the same as in the classical case $q=1$). 
\end{proposition} 

\begin{proof} $T_1\otimes T_2$ is a finite dimensional $\tuq$-module of type 
$1$, hence is the direct sum of modules from ${\cal{T}}$. 
The formal characters of $T_1$ and $T_2$ are the same as in the classical 
case $q=1$ by Proposition~\ref{prop:character-q-independent}. 
They determine the formal character of $T_1\otimes T_2$, so it is again the 
same as in the case $q=1$. 
So the assertion on the multiplicities follows by 
Proposition~\ref{prop:structure=character}. 
\end{proof}  

In the special case when $T_2$ is the vector representation of $\tuq$ 
(the irreducible representation with highest weight $(1,0,\ldots,0)$, 
the above result is proved in Proposition 4.2 (1) of \cite{h} 
(see also 8.6.2 Proposition 36 in \cite{ks}) 
by different 
methods. 

Finally note that in the odd dimensional case, 
the full orthogonal group $O(2l+1)$ is generated over $SO(2l+1)$ 
by the central element $-I$, which acts as a scalar $+1$ or $-1$ 
in any irreducible representation of $O(2l+1)$. 
Therefore the algebra $\cuu(O(2l+1))$ is defined as the 
tensor product $\mc[\chi]\otimes U_{q^{1/2}}({\mathrm{so}}_{2l+1})$, 
where $\chi$ here is just an abstract generator of the two-element group, 
and $\mc[\chi]$ is the corresponding group algebra. 
Then an irreducible $U_{q^{1/2}}({\mathrm{so}}_{2l+1})$-representation  
has always two extensions to an $\cuu(O(2l+1))$-representation on the same 
underlying space: the element $\chi$ acts as a scalar $+1$ or $-1$. 
The analogue of Proposition~\ref{prop:tensorproduct-multiplicities} 
holds obviously in this case. 


\section{Appendix B}\label{appendix-b} 

Here we sketch a proof of Theorem~\ref{thm:coc-gen-ogq} in the classical 
case $q=1$. 

When $G=SO(2l+1)$ or when $G$ is simple and simply connected, 
$\ogg\coc$ is a polynomial algebra generated by the characters of the 
fundamental representations, see for example \cite{s}. 
For $SL(l+1)$ or $SO(2l+1)$ the fundamental representations 
are the first $l$ exterior powers of the defining representation, hence we 
have (i) and (iv). The $r$th exterior power of the defining representation of 
$Sp(2l)$ for $r=1,\ldots,l$ is the direct sum of the $r$th fundamental 
representation and some copies of the fundamental representations with 
strictly lower index, see section 5.1.3 in \cite{gw}. 
Therefore $\sigma_1,\ldots,\sigma_l$ is another 
generating system of $\coo(Sp(2l))\coc$, and we get (ii). 
Since $O(2l+1)\cong SO(2l+1)\times\mz_2$, and $\omega_{2l+1}$ 
is trivial on $SO(2l+1)$ whereas it gives the non-trivial 
irreducible representation on $\mz_2$, the statement (iii) immediately 
follows from (iv). 

(v) Note that $G$ acts on itself by conjugation, and $\ogg\coc$ is the 
corresponding algebra of polynomial invariants. 
Identify $O(2l)$ with the subset of the space $M(2l,\mc)$ of 
$(2l\times 2l)$ 
matrices consisting of matrices $A$ with $AA^T=I$. The group $O(2l)$ 
acts on $M(2l,\mc)$ by conjugation, and the corresponding algebra of 
polynomial invariants is generated by the functions 
$A\mapsto \Tr(f(A,A^T))$ as $f$ ranges over the possible monomials in $A$ and 
$A^T$, see \cite{sib} or \cite{p}. 
Using $A^T=A^{-1}$ for $A\in O(2l)$, 
we get that the algebra $\ogg\coc$ is generated by 
the functions $A\mapsto \Tr(A^d)$, $d=1,\ldots,2l$ 
(the upper bound on $d$ comes from the Cayley-Hamilton identity). 
Note that $\sigma_i(A)$ is the $i$th characteristic coefficient of the matrix 
$A$, hence $\sigma_1,\ldots,\sigma_{2l}$ also generate $\ogg\coc$. 
For $r=1,\ldots,l$ we have the well known isomorphisms 
$\bigwedge^r\mc^{2l}\otimes \bigwedge^{2l}\mc^{2l}\cong
\bigwedge^{2l-r}\mc^{2l}$ of $O(2l)$-representations (see for example 
Exercise 6 in section 5.1.8 of \cite{gw}). 
This implies $\sigma_{2l-r}=\sigma_r\sigma_{2l}$, $r=1,\ldots,l$, 
hence $\ogg\coc$ is generated by $\sigma_1,\ldots,\sigma_l,\sigma_{2l}$, 
and the relations $\sigma_{2l}^2=1$ and $\sigma_l\sigma_{2l}=\sigma_l$ hold. 
We need to show that there are no further relations among these generators. 
Realize now $O(2l)$ as the group of invertible matrices 
$\{A\mid JA=(A^T)^{-1}J\}$, 
where $J=\left(\begin{matrix}0&1\\1&0\end{matrix}\right)
\oplus\cdots\oplus
\left(\begin{matrix}0&1\\1&0\end{matrix}\right)$, 
and restrict the functions in $\ogg$ to the subset 
$Y\sqcup Z$, where $Y$ consists of the diagonal matrices 
${\rm{diag}}(z_1,z_1^{-1},\ldots,z_l,z_l^{-1})$ 
and $Z$ consists of the matrices 
$\left(\begin{matrix}0&z_1\\z_1^{-1}&0\end{matrix}\right)
\oplus{\rm{diag}}(z_2,z_2^{-1},\ldots,z_l,z_l^{-1})$ with 
$z_i\in\mc^{\times}$.  
Now suppose that 
$f(\sigma_1,\ldots,\sigma_l)+\sigma_{2l}g(\sigma_1,\ldots,\sigma_{l-1})=0$ 
holds in $\ogg$. 
Clearly the restrictions of $\sigma_1,\ldots,\sigma_l$ to $Y$ are 
algebraically independent. 
So restricting the above relation to $Y$ we get that 
$f(t_1,\ldots,t_l)=-g(t_1,\ldots,t_{l-1})$ as polynomials in $t_1,\ldots,t_l$, 
so the above relation is 
$(1-\sigma_{2l})g(\sigma_1,\ldots,\sigma_{l-1})=0$. 
Now restricting this relation to $Z$ one sees that $g(t_1,\ldots,t_{l-1})$ is 
the zero polynomial. 

(vi) The highest weights of the representations corresponding to 
$\sigma_1,\ldots,\sigma_{l-1},\sigma_{l,0},\sigma_{l,1}$ 
generate the semigroup of the highest weights of all 
irreducible representations of $SO(2l)$, see for example pages 102 and 234 in 
\cite{gw}. 
Using the usual partial ordering on the weight semigroup, 
an inductive argument shows that the trace of an arbitrary irreducible 
$SO(2l)$-representation can be expressed as a polynomial of 
$\sigma_1,\ldots,\sigma_{l-1},\sigma_{l,0},\sigma_{l,1}$. 
The elements $\sigma_1,\ldots,\sigma_l$ are 
algebraically independent by the 
same argument as in (v). 
The full orthogonal group $O(2l)$ acts on $SO(2l)$ by conjugation, 
and this induces an 
action on $\coo(SO(2l))$. For $\chi\in O(2l)\setminus SO(2l)$ 
we have $\chi (\sigma_{l,0}-\sigma_{l,1})=-(\sigma_{l,0}-\sigma_{l,1})$, 
because the automorphism of $SO(2l)$ induced by $\chi$ interchanges the 
representations $\omega_{l,0}$ and $\omega_{l,1}$. 
Consequently, 
$(\sigma_{l,0}-\sigma_{l,1})^2$ is an $O(2l)$-invariant in $\coo(SO(2l))$, 
hence it is a polynomial of $\sigma_1,\ldots,\sigma_l$ by (v). 
The $\mc[\sigma_1,\ldots,\sigma_l]$-module generated by 
$1$ and $\sigma_{l,0}-\sigma_{l,1}$ is free (of rank two), 
because $\coo(SO(2l))$ is a domain, and the elements of 
$(\sigma_{l,0}-\sigma_{l,1})\mc[\sigma_1,\ldots,\sigma_l]$ are 
not invariant with respect to the action of 
the full orthogonal group, 
whereas $\sigma_1,\ldots,\sigma_l$ are $O(2l)$-invariants.  

The $SO(2l)$-invariant $\sigma_{l,0}-\sigma_{l,1}$ and the relation 
$(\sigma_{l,0}-\sigma_{l,1})^2=h(\sigma_1,\ldots,\sigma_l)$ 
can be seen more explicitly as follows. 
Think of $SO(2l)$ as the set of determinant $1$ matrices $A$ with 
$JA=(A^T)^{-1}J$. 
It is not difficult to check that up to sign, $\sigma_{l,0}-\sigma_{l,1}$ 
is the function mapping $A\in SO(2l)$ to the Pfaffian 
$\pf(JA-A^TJ)$ of the skew symmetric $(2l\times 2l)$ matrix 
$JA-A^TJ$. 
(For the definition and basic properties of the Pfaffian see 
the appendix of \cite{gw}; 
the $SO(2l)$-invariant $A\mapsto\pf(JA-A^TJ)$ appears in \cite{atz}.) 
Indeed, both $\sigma_{l,0}-\sigma_{l,1}$ and 
$A\mapsto\pf(JA-A^TJ)$ span a $1$-dimensional $O(2l)$-invariant 
subspace in $\coo(SO(2l))$ on which $O(2l)$ acts by the determinant 
representation. Both of them are contained in the space of matrix elements of 
the $l$th tensor power of the defining representation of $SO(2l)$, 
and in this $O(2l)$-invariant subspace of $\coo(SO(2l))$ 
the determinant representation of 
$O(2l)$ occurs with multiplicity one. So 
$\sigma_{l,0}-\sigma_{l,1}$ and $A\mapsto\pf(JA-A^TJ)$ are non-zero 
scalar multiples of each other. Restricting them to the maximal torus of 
$SO(2l)$ one can check that in fact they coincide (up to sign). 
For $A\in SO(2l)$ we have 
\[\pf^2(JA-A^TJ)=\det(JA-JA^{-1})=(-1)^l\det(A+I)\det(A-I).\]  
Therefore the relation 
\[(\sigma_{l,0}-\sigma_{l,1})^2
=(-1)^l
\left(\sigma_l+2\sum_{i=0}^{l-1}\sigma_i\right)
\left((-1)^l\sigma_l+2\sum_{i=0}^{l-1}(-1)^i\sigma_i\right)
\]
holds. 


\section{Appendix C}\label{appendix-c}

Here we deduce the assertion of 
Theorem~\ref{thm:coc-gen-frt} in the classical case $q=1$, for 
$G=O(N)$ or $Sp(N)$. 
Recall that $\clfrt$ is the coordinate ring $\coo(\mmm)$ of the 
Zariski closure $\mmm$ of the cone $\mc G$. 
The group $G$ acts on $M(N,\mc)$ by conjugation, 
and $\mmm$ is a $G$-stable subvariety in $M(N,\mc)$. 
We claim that $\clfrt\coc$ coincides with the algebra $\coo(\mmm)^G$ 
of $G$-invariants.  
This folows from the well known fact that for any affine algebraic group $H$, 
the algebra $\coo(H)\coc$ coincides with the algebra of adjoint invariants 
$\coo(H)^H$. Applying this for $H=\mc^{\times}G$, and observing 
that the conjugation action of $\mc^{\times}$ on $M(N,\mc)$ is trivial, 
we obtain that 
\[\clfrt\coc=\clfrt\cap\coo(\mc^{\times}G)\coc
=\coo(\mmm)\cap\coo(\mc^{\times}G)^{\mc^{\times}G}
=\coo(\mmm)\cap\coo(\mc^{\times}G)^G=\coo(\mmm)^G.\] 
Consider the natural surjection $\coo(M(N,\mc))^G\to\coo(\mmm)^G$. 
Generators of $\coo(M(N,\mc))^G$ are known from \cite{sib}, \cite{p}, 
these are the functions 
\[A\mapsto \Tr(A^{i_1}(A^*)^{j_1}\cdots A^{i_s}(A^*)^{j_s}),\]
where $A^*$ denotes the adjoint of $A\in M(N,\mc)$ with respect to the 
invariant bilinear form determining $G$.  
By definition of $\mmm$, if $A\in \mmm$, then 
$AA^*=A^*A$ equals the scalar  matrix $\ddd(A)I$. 
Note that $\Tr(A^*)=\Tr(A)$. Therefore, for $A\in\mmm$ we have 
\[\Tr(A^{i_1}(A^*)^{j_1}\cdots A^{i_s}(A^*)^{j_s})=\ddd(A)^k\Tr(A^d),\] 
where $k={\mathrm{min}}\{i_1+\cdots+i_s,j_1+\cdots+j_s\}$, 
and 
$d=|i_1+\cdots+i_s-j_1-\cdots-j_s|$. 
Taking into account the Cayley-Hamilton Theorem we get that $\coo(\mmm)^G$ 
is generated by the functions 
$A\mapsto \ddd(A)$, $A\mapsto\Tr(A^j)$, $j=1,\ldots,N$. 
By the Newton formulae the elements $\ddd,\tisi_1,\ldots,\tisi_N$ generate 
the same algebra. 

Next we determine the relations among the above generators. 
Identify $\clfrt$ with its image under the map 
$\iota:\clfrt\to\ogg\otimes\mc[z]$ 
from Proposition~\ref{prop:frt-injection}. 
Thus 
$\tisi_j=\sigma_jz^j$ and $\ddd=z^2$ 
(we supress the $\otimes$ sign from the notation). 
Since $\sigma_1,\ldots,\sigma_l$ are algebraically independent 
in $\ogg$ (see Theorem~\ref{thm:coc-gen-ogq}), the elements 
$\sigma_1z,\ldots,\sigma_lz,z^2$ are algebraically independent in $\clfrt$. 

When $G=Sp(N)$ $(N=2l)$, we have $\sigma_N=1$ and 
$\sigma_{N-i}=\sigma_i$ for $i=1,\ldots,l$ 
(this follows from the well known $G$-module isomorphism 
$\bigwedge^i\mc^N\otimes \bigwedge^N\mc^N=\bigwedge^{N-i}\mc^N$, 
and the fact that the 
$N$th exterior power of $\mc^N$ is the trivial $Sp(N)$-module).  
Thus  
\[\tisi_{N-i}=\sigma_{N-i}z^{N-i}=\sigma_iz^i(z^2)^{l-i}
=\tisi_i\ddd^{l-i}\] 
for $i=0,\ldots,l-1$. 
So $\caa(Sp(N))\coc$ is generated by $\ddd,\tisi_1,\ldots,\tisi_l$. 

Finally, assume $G=O(N)$. In $\coo(O(N))$ the relations 
$\sigma_N^2=1$ and 
$\sigma_i\sigma_N=\sigma_{N-i}$ for $i=1,\ldots,l$ hold, see 
Theorem~\ref{thm:coc-gen-ogq}. 
Therefore in $\clfrt$ we have 
\[\tisi_{N-i}\tisi_{N-j}=\sigma_{N-i}z^{N-i}\sigma_{N-j}z^{N-j}
=(\sigma_N)^2\sigma_iz^i\sigma_jz^j z^{2(N-i-j)}=\tisi_i\tisi_j\ddd^{N-i-j}\] 
for $0\leq i\leq j\leq l$, and 
\[\tisi_i\tisi_{N-j}\ddd^{j-i}=
\sigma_iz^i\sigma_{N-j}z^{N-j}z^{2(j-i)}
=\sigma_i\sigma_j\sigma_Nz^{N-i+j}
=\sigma_{N-i}z^{N-i}\sigma_jz^j
=\tisi_{N-i}\tisi_j\] 
for $0\leq i<j\leq l$. 
It is an elementary exercise to show that 
modulo these relations an arbitrary monomial of 
$\ddd,\tisi_1,\ldots,\tisi_N$ can be rewritten into a monomial contained in 
$B(N)$: using the relations of the first type we can get rid of those 
products of the generators which contain at least two factors from 
$\{\tisi_{l+1},\ldots,\tisi_N\}$. 
In the case $N=2l$, by the relation 
$\tisi_{N-i}\tisi_l=\tisi_i\tisi_l\ddd^{l-i}$ 
(the special case $j=l$ 
of the second type relations) 
we eliminate 
the products which 
contain $\tisi_l$ and a factor from $\{\tisi_{l+1},\ldots,\tisi_N\}$. 
By the relations $\tisi_{N-j}\ddd^j=\tisi_j\tisi_N$ (the special case $i=0$ 
of the second type relations) 
we get rid of the products which contain the 
subword $\tisi_{N-j}\ddd^j$ for $j=1,\ldots,l$ if $N=2l+1$ 
and for $j=1,\ldots,l-1$ if $N=2l$. 
Take a product of the generators which is not ruled out by the above 
reductions, and which is not contained in $B(N)$. 
Then it must contain a subword 
$\tisi_{N-j}\tisi_i\ddd^{j-i}$ with $1\leq i<j\leq l$ 
(respectively, $1\leq i<j\leq l-1$) when $N=2l+1$ (respectively, $N=2l$). 
Replace this subword by $\tisi_j\tisi_{N-i}$, using the second type relations. 
In this way we increase the index of the unique factor of this monomial from 
the set $\{\tisi_{l+1},\ldots,\tisi_{N-1}\}$. After finitely many such steps 
we end up in $B(N)$ or with a monomial eliminated already.  
So we have proved that the elements in $B(N)$ span $\caa(O(N))\coc$. 
We know from 
Theorem~\ref{thm:coc-gen-ogq} (iii) and (v) 
that $\sigma_1^{i_1}\cdots\sigma_l^{i_l}$, 
$\sigma_{2l+1}\sigma_1^{j_1}\cdots\sigma_l^{j_l}$ 
($i_s,j_s\in\mn_0$) are linearly independent in $\coo(O(2l+1))$, 
and 
$\sigma_1^{i_1}\cdots\sigma_l^{i_l}$,  
$\sigma_{2l}\sigma_1^{j_1}\cdots\sigma_{l-1}^{j_{l-1}}$ 
($i_s,j_s\in\mn_0$) are linearly independent in $\coo(O(2l))$. 
Using again the embedding $\iota$ we easily get that the elements of $B(N)$ 
are linearly independent in $\caa(O(N))$. 
Finally, the fact that $B(N)$ is a basis of $\caa(O(N))\coc$ 
implies that the set of relations used to rewrite arbitrary 
products of the generators 
as linear combinations of elements of $B(N)$ 
is complete: namely, the ideal of relations among the generators 
$\ddd,\tisi_1,\ldots,\tisi_N$ is generated by the relations given in the 
statement of Theorem~\ref{thm:coc-gen-frt}. 


\section{Appendix D}\label{appendix-d} 

Here we give a direct proof of the fact that for $G=O(N)$ or $Sp(N)$, 
the algebra $\clfrt$ defined in terms of generators and relations 
\eqref{eq:rel-clfrt} in Section~\ref{sec:coc-frt}, coincides with the 
coordinate ring of the Zariski closure $\mmm$ of the cone $\mc G$, 
where $G$ is embedded into $M(N,\mc)$ in the usual way; 
that is, $G=\{A\in M(N,\mc)\mid A\ccc^{-1} A^T\ccc=I\}$, 
where $\ccc$ is the matrix of a symmetric (respectively, skew-symmetric)  
non-degenerate bilinear form (the matrix $\ccc=\ccc(1)$ is specified in 
Section~\ref{sec:coc-frt}). 
In other words, we claim that the vanishing ideal $I(\mmm)$ of $\mmm$ 
in $\coo(M(N,\mc))=\mc[u^i_j\mid i,j=1,\ldots,N]$ is generated 
by the entries of 
$\kkk(1)\uu_1\uu_2-\uu_1\uu_2\kkk(1)$ (notation explained in 
Section~\ref{sec:coc-frt}). 
Write $B$ for the set of entries of this matrix, and write 
$\langle B\rangle$ for the ideal generated by these homogeneous 
quadratic elements. 
One sees directly from the definition of $\kkk(1)$ that 
$\langle B\rangle\subseteq I(\mmm)$, see for example the 
proof of 9.3.1 Lemma 12 in \cite{ks}. 
Since $\mc\mmm=\mmm$, the ideal $I(\mmm)$ is homogeneous. 
Take an arbitrary $f\in I(\mmm)$. Our aim is to show that $f$ is contained 
in $\langle B\rangle$. 
We may assume that $f$ is homogeneous of degree $d$. 
Clearly $f\in I(G)$, since $G\subset\mmm$. 
Now Theorems (5.2 C) and (6.3 B) of \cite{w} 
assert that $B$ and $\ddd-1$ generate $I(G)$ in a nice way; 
that is, there are elements $f_b,h\in\mc[u^i_j]$ $(b\in B)$, such that 
\begin{equation}\label{eq:ap-d-1}
f=(\ddd-1)h+\sum_{b\in B}bf_b,
\end{equation}
moreover, $\deg(f_b)\leq d-2$ and $\deg(h)\leq d-2$. 
We may assume that $h$ has the minimal possible number of non-zero 
homogeneous components. Suppose that $h\neq 0$. 
Write $h=\widehat{h}+\widetilde{h}$, where $\widetilde{h}$ is the minimum degree 
homogeneous component of $h$. Then 
\[(\ddd-1)h=-\widetilde{h}+\mbox{ higher degree terms.}\] 
Since $\deg(\widetilde{h})<d=\deg(f)$, it follows 
from \eqref{eq:ap-d-1} that $-\widetilde{h}$ is killed by the appropriate 
homogeneous component 
of $\sum_{b\in B}bf_b$, 
hence 
$\widetilde{h}=\sum_{b\in B}bh_b$ 
for some $h_b$, with $\deg(h_b)\leq d-4$. 
Thus we have 
\begin{equation}\label{eq:ap-d-2} 
f=(\ddd-1)\widehat{h}+
\sum_{b\in B}b(f_b+h_b(\ddd-1))
\end{equation} 
Note that in \eqref{eq:ap-d-2} we have 
$\deg(f_b+h_b(\ddd-1))\leq d-2$, 
and $\widehat{h}$ has fewer non-zero homogeneous components than 
$h$ in \eqref{eq:ap-d-1}. 
This contradiction implies that $h=0$ in \eqref{eq:ap-d-1}, 
so $f=\sum_{b\in B}bf_b$ is contained in $\langle B\rangle$.


\end{document}